\documentclass[11pt]{amsart}
\usepackage{amsmath,amssymb,amsfonts,url,mathptmx}
\numberwithin{equation}{section}

\newtheorem{thm}{Theorem}[section]
\newtheorem{lem}{Lemma}[section]
\newtheorem{rem}{Remark}[section]
\newtheorem{prop}{Proposition}[section]

\newcommand{\hdot}{^\text{\r{}}\hspace{-.33cm}H}

\begin{document}
\title[Liouville system]{On Liouville systems at critical parameters, Part 2: Multiple bubbles } 

\author{Hsin-Yuan Huang}
\address{Department of Applied Mathematics\\
National Chiao-Tung University\\
Hsinchu, Taiwan  } \email{hyhuang@math.nctu.edu.tw}

\author{Lei Zhang}\footnote{Lei Zhang is partially supported by a Simons Foundation Collaboration Grant}
\address{Department of Mathematics\\
        University of Florida\\
        358 Little Hall P.O.Box 118105\\
        Gainesville FL 32611-8105}
\email{leizhang@ufl.edu}

\date{\today}

\begin{abstract}
In this paper, we continue to consider the generalized Liouville system:
$$
\Delta_g u_i+\sum_{j=1}^n a_{ij}\rho_j\left(\frac{h_j e^{u_j}}{\int h_j e^{u_j}}- {1} \right)=0\quad\text{in \,}M,\quad i\in I=\{1,\cdots,n\},
$$
where $(M,g)$ is a Riemann surface  $M$ with volume $1$, $h_1,..,h_n$ are positive smooth functions and $\rho_j\in \mathbb R^+$($j\in I$). In previous works Lin-Zhang identified a family of hyper-surfaces $\Gamma_N$ and proved a priori estimates for $\rho=(\rho_1,..,\rho_n)$ in areas separated by $\Gamma_N$. Later Lin-Zhang also calculated the leading term of $\rho^k-\rho$ where $\rho\in \Gamma_1$ is the limit of $\rho^k$ on $\Gamma_1$ and $\rho^k$ is the parameter of a bubbling sequence. This leading term is particularly important for applications but it is very hard to be identified if $\rho^k$ tends to a higher order hypersurface $\Gamma_N$ ($N>1$). Over the years numerous attempts have failed but in this article we overcome all the stumbling blocks and completely solve the problem under the most general context: We not only capture the leading terms of $\rho^k-\rho\in \Gamma_N$, but also reveal new robustness  relations of coefficient functions at different blowup points.
\end{abstract}

\maketitle

\section{Introduction}
Let $(M,g)$ be a compact Riemann surface with volume $1$, $h_1,...,h_n$ be positive $C^3$ functions on $M$, $\rho_1,..,\rho_n$ be nonnegative
constants. In this article we continue our study of the following Liouville system defined on $(M,g)$:
\begin{equation}\label{mainsys}
\Delta_g u_i+\sum_{j=1}^n\rho_j a_{ij} \left(\frac{h_je^{u_j}}{\int_Mh_je^{u_j}dV_g}-1\right)=0,\quad i\in I:=\{1,..,n\}
\end{equation}
where   $dV_g$ is the volume form, $\Delta_g$ is the Laplace-Beltrami operator $\Delta_g\le 0$.
When $n=1$, equation \eqref{mainsys} is the mean field equation of the Liouville type:
\begin{equation}\label{equfromsys}
\Delta_g u+\rho\left(\frac{h e^u}{\int_M h e^udV_g}- {1} \right)=0\quad\text{in \,}M
\end{equation}
when $a_{11}=1$. Therefore, the Liouville system \eqref{mainsys} is a natural extension of the classical Liouville equation,
which has been extensively studied for decades because of its profound connections with various fields in geometry and physics. Since the general form of Liouville systems includes many models from Biology,  Physics  and other disciplines of sciences, it is very desirable to study generical Liouville systems and derive common features.
Recently, the Liouville system
has drawn a lot of attention because it also arises from
the stationery solutions of  multi-species Patlak-Keller-Segel system\cite{wolansky3}  and self-dual condensate solutions of Ablian Chern-Simons  model with $N$ Higgs particles\cite{Wil,kimleelee} when certain parameter tends to zero. In particular,  these two examples exhibit the bubbling phenomenon.  The study of bubbling solutions represents an essential difficulty of Liouville system and it not only impacts the immediately related fields but also depends on the development of them. The readers may look into the following references for closely related discussions \cite{aly,barto1,barto3,bennet,caffarelli,lion1,lion2,ccl,chanillo2,chchlin,childress,debye,dunne,dziar,phys,huang-zhang,hong,jackiw,keller,kiess,kiessling2,nolasco2,rubinstein,spruck,wolansky1,
wolansky2,yang,zhangcmp,zhangccm}.\par

For system (\ref{mainsys}) the Sobolev spaces for solutions are
$$\hdot^{1,n}=\,\hdot^1(M)\times\cdots\times\,\hdot^1(M) $$
where
\begin{equation*}
\hdot^1(M)=\left\{
u\in L^2(M)\,\Big|\, |\nabla_g u|\in L^2(M)\text{ and }\int_M  u\,dV_g=0
\right\}.
\end{equation*}
For any $\rho=(\rho_1,\cdots, \rho_n)$, $\rho_i>0$, let $\varPhi_\rho$ be a nonlinear functional defined in \, $\hdot^{1,n}$ by
$$
\varPhi_\rho(u)=\frac{1}{2}\sum_{i=1}^n\sum_{j=1}^n a^{ij}\int_M \nabla_g u_i\cdot \nabla_g u_j dV_g- \sum_{j=1}^n\rho_j\log\int_M h_j e^{u_j}dV_g
$$
where $(a^{ij})_{n\times n}$ is the inverse of $A=(a_{ij})_{n\times n}$, $I$ is the set of indexes: $I=\{1,...,n\}$.
It is easy to see that equation \eqref{mainsys} is the Euler-Lagrangian equation of $\varPhi_\rho$.

In \cite{linzhang1,linzhang2}, Lin and the second  author   completed a degree counting program for (\ref{equfromsys}) under the following two assumptions on the
matrix A:

\begin{align*}
&(H1): \quad A\mbox{ is symmetric, nonnegative, irreducible and invertible.}\\
&(H2): \quad a^{ii}\leq 0,\,\, \forall i\in I,\quad  a^{ij}\geq 0 \,\, \forall i\neq j\in I,  \quad  \sum_{j=1}^na^{ij}\geq 0 \,\,\forall i\in I.
\end{align*}
Roughly speaking $(H1)$ is a rather standard assumption for Liouville systems, $(H2)$ says the interaction between equations has to be strong.
In \cite{linzhang2} Lin-Zhang identified a family of hypersurfaces
\begin{equation*}
\Gamma_N=
\left\{
\rho\,\big|\, \rho_i>0, i\in I;\,\,  \Lambda_I(\rho)=0,\quad \Lambda_J(\rho)>0,\,\, \forall \,\, \emptyset \neq J\subsetneq I,\,\,
\right\}.
\end{equation*}
where
$$\Lambda_I(\rho)=4\sum_{i=1}^n\frac{\rho_i}{2\pi N}-\sum_{i=1}^n\sum_{j=1}^na_{ij}\frac{\rho_i}{2\pi N}\frac{\rho_j}{2\pi N}$$ ($\Lambda_J(\rho)$ is understood similarly).
Furthermore they proved that if $\rho=(\rho_1,...,\rho_n)$ stays in a regions bounded by these critical hypersurfaces, a priori estimate holds. Based on the a priori estimate, Lin-Zhang proved the following degree counting formula which implies existence of solution if the degree is not zero:
\begin{equation}\label{dcount}
d_{\rho}=\left\{\begin{array}{ll}
1\quad \text{if }\rho \in \mathcal{O}_0\\
\frac{1}{N!}\bigg ((-\chi_M+1)...(-\chi_M+N) \bigg )\quad \text{if }\rho \in \mathcal{O}_N.
\end{array}
\right.
\end{equation}
where $\chi_M$ is the Euler characteristic of $M$, $\mathcal{O}_N$ is the region between $\Gamma_N$ and $\Gamma_{N+1}$. The jump of the degree from $\mathcal{O}_{N-1}$
to $\mathcal{O}_{N}$ is  contributed by the  blowup solutions to \eqref{mainsys}.\par

In this article we consider the case $\rho\to \Gamma_N$ from $\mathcal{O}_{N-1}$ or $\mathcal{O}_N$: suppose $\rho=(\rho_1,..,\rho_n)\in \Gamma_N$ is a limit point,
$\rho^k=(\rho_1^k,..,\rho_n^k)$ is a sequence of parameters corresponding to bubbling solutions, the aim of this article is to identity the leading term of $\rho^k-\rho$ as $k\to \infty$.  Since the normal vector at $\rho$ is proportional to
$$( \sum_{j=1}^na_{1j}\frac{\rho_j}{2\pi N}-2,..., \sum_{j=1}^na_{nj}\frac{\rho_j}{2\pi N}-2), $$
which has all components positive (implied by $\Lambda_J>0$ in the definition of $\Gamma_N$, see \cite{linzhang2}),  we assume
\begin{equation}\label{rho-com}
\frac{\rho_i^k-\rho_i}{\rho_j^k-\rho_j}\sim 1\quad \forall  \quad i\neq j\in I.
\end{equation}
Note that $A_k\sim B_k$ means $CA_k\le B_k\le C_1A_k$ for some $C,C_1>0$ independent of $k$. It is established in \cite{linzhang2} that when $\rho^k=(\rho_1^k,..,\rho_n^k)$ tends to $\rho\in\Gamma_N$, blowup solutions have exactly $N$ disjoint blowup points: $p_1,...,p_N$. In \cite{lin-zhang-jfa} Lin-Zhang derived the leading term when $N=1$. It is interesting to observe in \cite{lin-zhang-jfa} that there is one particular point $Q\in\Gamma_1$ such that if $\rho\to Q$ the leading term contains local curvature at $Q$ only, but if $\rho$ tends to any other point, the leading term is involved with global integration of the whole manifold.\par

The main purpose of this article is to extend the result in \cite{lin-zhang-jfa} to $\rho^k\to \Gamma_N$ when $N>1$.
Among   other things, we obtain the leading term of $\Lambda_{I}(\rho^k)$ as $\rho^k\to \Gamma_N$   which gives us the sufficient conditions of a uniform bound of solutions as  $\rho^k\to \Gamma_N$. Let $p_1$,...,$p_N$ be $N$ disjoint blowup points, which means for each $p_t$, there exist $p_t^k\to p_t$ such that $  \max_{i\in I}u_i^k(p_t^k)\to \infty$, $t\in\{1,\cdots,N\}$.  For $h_i^k$ we assume that they are uniformly bounded by positive constants:
\begin{equation}\label{assum-h}
\frac 1C\le h_i^k(x)\le C, \quad \|h_i^k\|_{C^3(M)}\le C
\end{equation}
for all $i$ and a $C>0$ independent of $k$.
Throughout the paper we use $u^k=(u_1^k,...,u_n^k)$ to denote blowup solutions and $M_{k,t}$ to denote the magnitude of $u^k$ near the blowup point $p_t$ and use $\epsilon_{k,t}=e^{-\frac 12M_{k,t}}$ to measure the errors:
\begin{equation}\label{Mk}
M_{k,t}=\max_{i\in I}\max_{x\in B(p_t,\delta_0)} \{u_i^k(x)-\log \int_M h_i^k e^{u_i^k}dV_g,\},\quad \epsilon_{k,t}=e^{-\frac 12 M_{k,t}},
\end{equation}
Here we require $\delta_0$ to be small enough so that $B(p_{t_1},\delta_0)\cap B(p_{t_2},\delta_0)=\emptyset $ for all $t_1\neq t_2$.  Let $p_t^k$ be the point where
 $$\max_{ i\in I}\max_{x\in B(p_t, \delta_0)} u_i^k(x)-\log \int_M h_i^k e^{u_i^k} dV_g $$
is taken in
$B(p_t, \delta_0)$.  Under the assumptions $(H1)$ and $(H2)$, we have a full blow-up picture in all balls (see \cite{linzhang2}).\par
To understand more precise information for the blowup phenomenon to \eqref{mainsys}, we shall study the convergence rate of $\Lambda_I(\rho^k)$ as $\rho^k\to \rho$
in terms of the magnitude of $u^k$.
The following  two fundamental questions will be answered in this paper:\\
\\
\begin{itemize}
\item[(1)] {\it Are the  magnitude of $u^k$  at different blowup points comparable to each other?}\\
\item[(2)] {\it What is the convergence rate of    the  difference of the local masses $\sigma_{i,s}^k-\sigma_{i,t}^k$ where
$\sigma_{i,s}^k=\int_{B( p_{s},\delta_0)}e^{u_i^{k}}dV_g$ and $\sigma_{i,t}^k=\int_{B( p_{t},\delta_0)}e^{u_i^{k}}dV_g$?}\par
\end{itemize}

     The first main result is to answer the first question.
\begin{thm} \label{heightthm} Let $u^k\in \,\, \hdot^{1,n}(M)$ be a sequence of blowup solutions of (\ref{mainsys}). Suppose $(H1),(H2)$ holds for $A$, \eqref{rho-com} holds for $\rho^k$ and (\ref{assum-h}) holds for $h_i^k$. Then
\begin{equation}
|M_{k,s}-M_{k,t}|=O(1),\quad\text{ for }\,\,s\not=t,\,\, s,t\in\{1,\cdots, N\} .
\end{equation}
Here, $O(1)$ is independent of $k$.
\end{thm}
 Note that without knowing
$|M_{k,s}-M_{k,t}|=O(1) $ we do not even know if $O(\epsilon_{k,s})=O( \epsilon_{k,t})$. We thus can use
   $\epsilon_{k }=e^{-\frac 12 M_{k }}$ to measure the errors, where
   $$
M_{k }=\max_{i\in I}\max_{x\in M} \{u_i^k(x)-\log \int_M h_i^k e^{u_i^k}dV_g,\}.$$

Let
\begin{equation}\label{def-m}
m_i=\sum_{j=1}^n a_{ij}\frac{\rho_j}{2\pi N}, \quad i=1,..,n,\quad m=\min_{i\in I} m_i.
\end{equation}
Here we note that either $2<m<4$ or all $m_i=4$ for all $m_i$ (see (\ref{mle4eq4})). Now we define a special point $Q_N=(q_1,...,q_n)$ on $\Gamma_N$, which satisfies
$$ \sum_{j=1}^na_{ij}q_j=8\pi N \quad \forall i\in I. $$
The second result is showing the tightness of the local masses.
\begin{thm}\label{localmass} Under the same assumptions in Theorem \ref{heightthm}.
\begin{itemize}
\item[(1)] If
$\rho^k\to \rho\in \Gamma_N (\rho\neq Q_N)$ from $\mathcal{O}_{N-1}$ or $\mathcal{O}_N$, then
\begin{equation}
|\sigma_{i,s}^{k}-\sigma_{i,t}^k|=O(\epsilon_{k }^{m-2}),\quad \text{ for }s\not=t,\,\,s,t\in\{1,\cdots,N\},\,\, i\in I.
\end{equation}
\item[(2)] If
$\rho^k\to  Q_N  $ from $\mathcal{O}_{N-1}$ or $\mathcal{O}_N$, then
\begin{equation}
|\sigma_{i,s}^{k}-\sigma_{i,t}^k|=O(\epsilon_{k }^{ 2}\log\frac{1}{\epsilon_{k }}),\quad \text{ for }s\not=t,\,\,s,t\in\{1,\cdots,N\},\,\, i\in I.
\end{equation}

\end{itemize}
\end{thm}
We remark here that the techniques developed in the proof of Theorems \ref{heightthm} and \ref{localmass} also play the key roles in the  study of the local uniqueness of the bubbling solutions in Liouville systems\cite{HZ2021}. For the sake of contradiction,   one also need to compare two sequence of the bubbling solutions at the blowup points with the same limit $\rho$. \par

The leading terms of $\Lambda_I(\rho^k)$ are different in  two cases. Before we state the results, we give some notations:\\

We define $N$ open sets $\Omega_{t,\delta_0}$ such that they are mutually disjoint, each of them contains a bubbling disk and their union is $M$:
\begin{equation}\label{set-decom}
B(p_t^k,\delta_0)\subset \Omega_{t,\delta_0}, \quad \cup_{t=1}^N \overline{\Omega_{t,\delta_0}}=M,
\quad \Omega_{t,\delta_0}\cap \Omega_{s,\delta_0}=\emptyset, \,\, \forall t\neq s.
\end{equation}
Let
$$I_1=\{i\in I;\quad \lim_{k\to \infty}m_i^k=m.\quad \}. $$
and $G$ be the Green's function defined by
$$
-\Delta_g G(x,\cdot)=\delta_p-1,\quad
\int_MG(x,\eta)dV_g(\eta)=0,
$$
and $\gamma$ is the regular part of the Green's function. Note that in local coordinates of a point, say $\eta$, $G$ is of the form
$$G(x,\eta)=-\frac 1{2\pi} \log |x-\eta| +\gamma(x,\eta). $$
We also define \begin{equation}\label{gstar}
G^*(p_t^k,p_s^k)=\left\{\begin{array}{ll}
\gamma(p_t^k,p_t^k),\quad s=t, \\
G(p_t^k,p_s^k),\quad s\neq t.
\end{array}
\right. \qquad s,t =1,...,N.
\end{equation}

 The third result is   the leading terms of $\Lambda_I(\rho^k)$ of the first case.

\begin{thm} \label{gammaNpneq} Under the same assumptions in Theorem \ref{heightthm}.  If
$\rho^k\to \rho\in \Gamma_N (\rho\neq Q_N)$ from $\mathcal{O}_{N-1}$ or $\mathcal{O}_N$, then
\begin{equation}\label{11july13e9}
\Lambda_I(\rho^k)
=(D+o(1))\frac{\epsilon_k^{m-2}}{N}.
\end{equation}
 Here,
The quantity $D$ is defined as follows
\begin{equation}
   D
 = \sum_{i\in I_1}e^{D_i-\alpha_i}\sum_{t=1}^Nc_t\lim_{\delta_0\to 0}\bigg (\delta_0^{2-m}-\frac{(m-2)}{2\pi}\int_{\hat{\Omega}_{t,\delta_0}}
(\frac{h_i^k(x)}{h_i^k(p_t^k)}e^{2\pi m\sum_{l=1}^N(G(x,p_l^k)-G^*(p_t^k,p_l^k))}dV_g\bigg ).
\end{equation}
where  $\hat{\Omega}_{t,\delta_0}= \Omega _{t,\delta_0}\setminus B(p_t^k,\delta_0)$,
$c_t=\frac{h_i^k(p_t^k)e^{2\pi m\sum_{l=1}^NG^*(p_t^k,p_l^k)}}{h_i^k(p_1^k)e^{2\pi m\sum_{l=1}^N G^*(p_1^k,p_l^k)}}$ and $D_i$, $\alpha_i$ are constants defined in \eqref{U-global}.
\end{thm}
We remark here that $\Lambda_I(\rho^k)>0$ for $\rho^k\in \mathcal{O}_{N-1}$ and $\Lambda_I(\rho^k)<0$ for $\rho^k\in \mathcal{O}_{N}$. Thus, if $D\not=0$, then the blowup solutions with $\rho^k$   occur   as $\rho^k\to \rho\in\Gamma_N(\rho\not= Q_N)$ only from one side of $\Gamma_N$.  Furthermore, it  yields a uniform bound of solutions as $\rho^k$ converges
to $\Gamma_N(\rho\not= Q_N)$ from $\mathcal{O}_{N-1}$    provided that $D<0$.  It is easy to see that for fixed $k$ the following limit exists:
\begin{equation}\label{leadterm1}
\lim_{\delta_0\to 0}\bigg (\delta_0^{2-m}-\frac {m-2}{2\pi}\int_{\hat{\Omega}_{t,\delta_0}}\frac{h_i^k(x)}{h_i^k(p_t^k)}
e^{2\pi m(\sum_{l=1}^NG(x,p_l^k)-G^*(p_t^k,p_l^k))}dV_g \bigg )
\end{equation}
because the leading terms from the  Green's function is $-\frac 1{2\pi}\log |x-p_t^k|$. The study of the sign of $D$ is another interesting fundamental question. However, it is  out of the scope of of the present  article. We will come back to this issue in the further study.\par

The fourth major result is concerned with the leading term of
$\Lambda_I(\rho^k)$ when $\rho^k\to Q_N$ and the major difference is that the leading term only depends on the curvature and coefficient functions at blowup points: $p_1^k$, ..., $p_N^k$:
\begin{thm}\label{gammaNpeq} Under the same assumptions in Theorem \ref{heightthm}. If   $\rho^k\to Q_N$ from $\mathcal{O}_{N-1}$ or $\mathcal{O}_N$, then
\begin{equation}\label{11july13e10}
\Lambda_I(\rho^k)
=-4\sum_{i=1}^n\sum_{t=1}^N b_{it}^k\epsilon_k^2\log \epsilon_k^{-1}+O(\epsilon_k^2)
\end{equation}
where

\begin{align}
b_{it}^k&=e^{D_i-\alpha_i}\bigg (\frac 14\frac{\Delta h_i^k(p_t^k)}{h_i^k(p_t^k)}-K(p_t^k)+4\pi N\\
&+4\pi\frac{\nabla h_i^k(p_t^k)}{h_i^k(p_t^k)}\cdot \sum_{l=1}^N\nabla_1G^*(p_t^k,p_l^k)
+16\pi^2|\sum_{l=1}^N\nabla_1G^*(p_t^k,p_l^k)|^2\bigg ). \nonumber
\end{align}
and $K$ is   the Gaussian curvature.
\end{thm}

As an application of the formula \eqref{11july13e10}, we obtain the uniform bound of solutions as $\rho^k$ converges
to $   Q_N $ from $\mathcal{O}_{N-1}$ provided $b_{it}^k>0$ for $t=1,\cdots,N.$

The fifth result is about the locations of the blowup points and the mutual relation of coefficient functions.

\begin{thm}\label{main-3}
If $\rho^k\to \rho\neq Q$ as in (\ref{rho-com}), then for $t=1,...,N$
\begin{equation}\label{11july13e7}
\sum_{i=1}^n \big (\nabla (\log h_i^k)(p_t^k)+2\pi m_i\sum_{s=1}^N \nabla_1 G^*(p_t^k,p_s^k)\big )=O(\epsilon_k^{m-2}), \,\, t=1,..,N
\end{equation}
where $\nabla_1$ means the differentiation with respect to the first component.
If $\rho^k\to Q$ as in (\ref{rho-com}),
\begin{equation}\label{11july13e8}
\sum_{i=1}^n \big (\nabla (\log h_i^k)(p_t^k)+8\pi \sum_{s=1}^N \nabla_1 G^*(p_t^k,p_s^k)\big )=O(\epsilon_k^{2}\log \epsilon_k^{-1}), \,\, t=1,..,N
\end{equation}
\end{thm}

Besides the location of blowup points, we also reveal new information about coefficient functions in the next theorem. For convenience we set
\begin{equation}\label{htk}
H_{i,t}^k:=\log h_i^k(p_t^k)+2\pi m_i\sum_{l=1}^NG^*(p_t^k,p_l^k).
\end{equation}

\begin{thm}\label{mut-com} Let $H_{i,t}^k$ be defined as in (\ref{htk}), then we have
\begin{equation}\label{mut-re}
H_{i,t}^k=H_{i,s}^k+O(\epsilon_k^{m-2}),\quad \forall i\in I, \quad \forall t\neq s, \quad \mbox{if} \quad  m<4,
\end{equation}
\begin{equation}\label{mut-re-2}
H_{i,t}^k=H_{i,s}^k+O(\epsilon_k^{2}\log \frac{1}{\epsilon_k}),\quad \forall i\in I,\quad t\neq s, \quad \mbox{if}\quad m=4.
\end{equation}
\end{thm}

The tightness of the coefficient functions is not seen in
the case of $N=1$. This show that the occurrence of multi bubbles forces the coefficient functions at the blowup points to be the same.   In construction of bubbling solutions one needs to know the precise information about bubbling interactions, exact location of blowup points, accurate vanishing rate of coefficient functions and specific leading terms in asymptotic expansions. All these have been covered in the main results of this article. Until now the construction of bubbling solutions for Liouville systems is still in the early stage of development, as the constructions so far are still restricted to single blowup point situations \cite{huang-cvpde}.\par

In all these main results the readers can see that sharp estimates are obtained for all the error terms. This is why we think these results will play as a central role in applications.  We expect the theorems of this article to serve as a benchmark for more sophisticated discussions in the near future.
 Also Gu-Zhang \cite{gu-zhang} completed a degree counting program for singular Liouville systems, the corresponding discussion for leading terms of approximating critical hyper-surfaces of singular Liouville systems is another exciting unconquered land to explore. Besides these immediate impacts to closely related fields, the idea of the proof in this article, the way to overcome major difficulties in bubble interaction could lead to major advance in Chern-Simons type equations (see \cite{huang-zhang}) as well.\par

   The main difficulty in the proof of the main theorems is on the interaction of bubbling solutions, which has a large to do with the nature of Liouville systems.
   For global solutions defined in $\mathbb R^2$, it is established in \cite{CSW, linzhang1} that total integrations of all components form a $n-1$ dimensional hyper-surface similar to $\Gamma_1$. This continuum of energy brings great difficulty for bubbling interaction: If we only have two bubbling disks, the energy in each disk is very close to $\Gamma_1$, but to identify the leading term of $\rho^k-\rho\in \Gamma_2$ one has to prove that they are both tending the \emph{same point} on $\Gamma_1$. It does not help to use the fact that they have the same limit, because the energy sequence may tend to its limit position very slowly. This
   difficulty does not exist for Toda systems, because the energy set for Toda systems is discrete. When we have only one bubble, this bubble interaction situation can be avoided ( see \cite{lin-zhang-jfa}). So the main contribution in this article is to prove that the bubbling solutions have almost the same energy in each bubbling disk.  The key idea is as follows: Around each blowup point, we first use an approximation theorem of Lin-Zhang \cite{lin-zhang-jfa} to have an initial expansion of the bubbling solution. The first term of this expansion is a sequence of global solutions. In this article we identify what the global sequence is around each bubbling point and compare them using a key idea of Lin-Zhang \cite{linzhang1} in their proof of the classification theorem for Liouville systems.
It turns out that after scaling, the global sequences are extremely close to one another. Then we further prove that the energy of the global sequence is not too much different from the bubbling solutions in each bubbling disk. All the error terms must be carefully identified in order to single out the leading terms in the main theorems.

The organization of this article is as follows: In section two we deploy the basic setting for all the topics in this article and invoke a few approximation results in previous works of Lin-Zhang. In section three we prove the closeness of bubbling solutions around different blowup points, in which  Theorems \ref{heightthm}, \ref{localmass} and \ref{mut-com} will be proved. The proof is set-up in two stages as the approximation becomes better in the second stage.
Finally in section four all other main theorems will be proved based on the precise estimates in section three.

\section{Approximation around a blowup point}
First we claim that we can assume $u^k=(u_1^k,...,u_n^k)$ to satisfy
\begin{equation}\label{1013e2}
 \int_Mh_i^ke^{u_i^k}dV_g=1,\quad Vol(M)=1.
 \end{equation}
 because otherwise we just consider
\begin{equation}\label{11jun25e1}
\Theta_i^k=u_i^k-\log \int_M h_i^ke^{u_i^k}dV_g,\quad i\in I.
\end{equation}
 Then we have
 \begin{equation}\label{1013e1}
 -\Delta_g \Theta_i^k=\sum_{j=1}^na_{ij}\rho_j^k(h_j^ke^{\Theta_j^k}-1).
 \end{equation}
 where $\Theta^k$ satisfies (\ref{1013e2}). Here we first set up preliminary discussions about the profile of $u^k$ near a blowup point.
 Suppose $p$ is a blowup point and in $B(p,\delta)$ there is only one blowup point of $u^k$. Let
 $$\tilde M_k=\max_{i\in I, x\in B(p,\delta)}u_i^k(x)+\log (\rho_i^k h_i^k(\tilde p_k)) \quad \mbox{and }\quad \tilde \epsilon_k=e^{-\frac 12 \tilde M_k},$$
 and $\tilde p_k$ be where $\tilde M_k$ is attained ( $\tilde p_k\to p$).
  Then the functions
 $$\tilde v_i^k(y)=u_i^k(\tilde p_k+\tilde \epsilon_ky)+2\log \tilde \epsilon_k $$
 converge in $C^2_{loc}(\mathbb R^2)$ to the limit function $v=(v_1,...,v_n)$ which is a global solution of

\begin{equation}\label{817e1}
\left\{\begin{array}{ll}-\Delta v_i=\sum_{j}a_{ij}e^{v_j},\quad \mathbb R^2,\quad i=1,..,n\\
\\
\int_{\mathbb R^2}e^{v_i}<\infty,\quad i=1,..,n,\quad   \max_{i\in I} v_i(0)=0.
\end{array}
\right.
\end{equation}
Here we note that it is established in \cite{linzhang2} that with assumptions (H1), (H2) all the bubbling solutions are fully bubbling: the limit must have $n$ equations and no component is lost in the limiting taking process. The classification of all global solutions of (\ref{817e1}) was completed in the work of Chipot-Shafrir-Wolansky \cite{CSW} and Lin-Zhang \cite{linzhang1}. All components of $v=(v_1,...,v_n)$ have one common point of symmetric symmetry. In this context, this common point is the origin.

To state more precise approximation results we write the equation in local coordinates around $\tilde p_k$.
In this coordinate $ds^2$ has the form $e^{\phi(y_{\tilde p_k})}(dy_1^2+dy_2^2)$ where
\begin{equation}\label{11july11e1}
|\nabla \phi(0)|=0,\quad \phi(0)=0,\quad \Delta_{y_{\tilde p_k}}\phi=-2Ke^{\phi}
\end{equation}
 where $K$ is the Gauss curvature. In local coordinates, (\ref{11jun25e1}) becomes
 \begin{equation}\label{11may19e3}
 -\Delta u_i^k=\sum_{j=1}^na_{ij}\rho_j^ke^{\phi}(h_j^ke^{u_j^k}-1),  \quad \mbox{in}\quad B(0,\delta).
 \end{equation}
 Let $f_i^k$ be defined by
 \begin{equation}\label{11may19e4}\left\{\begin{array}{ll}
 \Delta f_i^k=\sum_{j=1}^na_{ij}\rho_j^ke^{\phi}, \quad \mbox{in}\quad B(0,\delta).\\
 \\
 f_i^k(0)=0, \quad \nabla f_i^k(0)=0.
 \end{array}
 \right.
 \end{equation}
 Then we have
\begin{equation}\label{11may19e5}
-\Delta (u_i^k-f_i^k)=\sum_{j=1}^na_{ij}\rho_j^kh_j^ke^{u_j^k-f_j^k}e^{f_j^k}e^{\phi},\quad \mbox{in}\quad B(0,\delta).
\end{equation}
If we set
$$\tilde h_i^k(x)=\frac{h_i^k(x)}{h_i^k(\tilde p_k)}e^{\phi+f_i^k}, $$
we have $\tilde h_i^k(0)=1$,
\begin{equation}\label{grad-h}
\nabla (\log \tilde h_i^k)(0)=\nabla (\log h_i^k)(\tilde p_k).
\end{equation}
 \begin{equation}\label{11july7e3}
\Delta (\log \tilde h_i^k)(0)=\Delta (\log h_i^k)(\tilde p_k)-2K(\tilde p_k)+\sum_{j=1}^n a_{ij}\rho_j^k.
\end{equation}
Thus we set
\begin{equation}\label{util-k}
\tilde u_i^k=u_i^k+\log \rho_i^k+\log h_i^k(\tilde p_k)-f_i^k,
\end{equation}
and write the equation of $\tilde u_i^k$ as
\begin{equation}\label{tuik}
\Delta \tilde u_i^k+\sum_{j=1}^n a_{ij}\tilde h_j^k e^{\tilde u_j^k}=0,  \quad \mbox{ in}\quad B(0,\delta).
\end{equation}

Now we introduce $\phi_i^k$ to be a harmonic function defined by the oscillation of $\tilde u_i^k$ on $B(\tilde p_k,\delta)$:
\begin{equation}\label{816e5}
\left\{\begin{array}{ll}
-\Delta \phi_i^k=0,\quad\mbox{in}\quad B(0,\delta),\\
\\
\phi_i^k=\tilde u_i^k-\frac 1{2\pi \delta}\int_{\partial B(0,\delta)}\tilde u_i^k,\quad \mbox{on}\quad \partial B(0,\delta).
\end{array}
\right.
\end{equation}
Obviously $\phi_i^k(0)=0$ by the mean value theorem and $\phi_i^k$ is uniformly bounded on $B(0,\delta/2)$ because $u_i^k$ has finite oscillation away from blowup points.
It is a standard fact (see \cite{linzhang1,linzhang2}) that the location of $\max_{i\in I}\max_{ x\in B(0,\delta)}\tilde u_i^k(x)-\phi_i^k(x)$ is $O(\epsilon_k^2)$ (roughly speaking, the reason is $0$ is a non-degenerate maximum of $ u_i^k$ and $\phi_i^k(0)=0$).

Going back to the original coordinate system, we call the maximum point after perturbation $p_k$. Now we set
$$M_k=\max_{i\in I}\max_{  x\in B(0,\delta)} u_i^k(x)+\log (\rho_i^kh_i^k(p_k))-\phi_i^k(x), \quad \epsilon_k=e^{-\frac 12 M_k}, $$
and we let
 $V^k=(V_1^k,..,V_n^k)$ be the radial solutions of
\begin{equation}\label{88e2}
\left\{\begin{array}{ll}-\Delta V_i^k=\sum_{j=1}^na_{ij}e^{V_j^k}\quad\mbox{in}\quad \mathbb R^2,\quad i\in I,\\  \\
V_i^k(0)=u_i^k(p_k)+\log (\rho_i^k h_i^k(p_k))-\phi_i^k(p_k),\quad i\in I.
\end{array}
\right.
\end{equation}
Note that since $p_k=\tilde p_k+O(\epsilon_k^2)$, it is easy to obtain that the oscillation of $V_i^k$ on $\partial B(p_k,\delta)$ is $O(\epsilon_k^2)$. The sequence of function $V^k=(V_1^k,...,V_n^k)$, which agrees with $u_i^k(x)+\log (\rho_i^kh_i^k(p_k))-\phi_i^k(x)$ at $p_k$, gives the first term in the approximation of $u_i^k$ near $p$. To state more precise approximation terms, we use
$$v_i^k(y)=u_i^k(p_k+\epsilon_k y)+\log (\rho_i^k h_i^k(p_k))-\phi_i^k(\epsilon_k y)+2\log \epsilon_k, \quad |y|<\frac{\delta}2 \epsilon_k^{-1} $$
and the following rough approximation theorem is established in \cite{lin-zhang-jfa}:

\begin{rem} The notations $M_k$, $\epsilon_k$ are the same as those in the introduction. It is confusing at this moment, later in the multiple bubbling situation we
will use $M_{k,t}$ and $\epsilon_{k,t}$ to denote the maximum of bubbling solutions and decay rate in each bubbling disk $B(p_t^k,\delta_0)$. $M_k$ is the maximum of $M_{k,t}$. But our analysis will show that we can replace $M_k$ by $M_{k,t}$, $\epsilon_k$ by $\epsilon_{k,t}$ for any $t$ and the nature of the proof does not change.
\end{rem}

Before citing the approximation theorems in \cite{lin-zhang-jfa} we mention one simple fact implied by the Pohozaev identity. Let
$\sigma_i^k=\frac 1{2\pi}\int_{B(p_k,\delta)}h_i^ke^{u_i^k}$ and $m_i^k= \sum_{j=1}^na_{ij}\sigma_j^k$. Let $\sigma_i=\lim_{k\to \infty}\sigma_i^k$ and
$m_i=\lim_{k\to \infty}m_i^k$. As usual we set $m=\min\{m_1,...,m_n\}$. Then it is established in \cite{linzhang1} that each $m_i>2$ and
$$\sum_{i=1}^n\sigma_i(m_i-4)=0.$$
Since each $\sigma_i>0$ it is easy to see that either
\begin{equation}\label{mle4eq4}
m<4, \quad \mbox{or} \quad m_i=4 \,\, \forall i\in I.
\end{equation}

The first approximation theorem established in \cite{lin-zhang-jfa} is a rough one that does not distinguish $m<4$ or $m=4$.
 \begin{thm}\label{cruexp} Given $\delta>0$, there exist $C(\delta)>0$, $k_0(\delta)>1$ such that for $|y|\le \frac{\delta}2 \epsilon_k^{-1}$ and $|\alpha |=0,1$, the following holds for all $k\ge k_0$
\begin{equation}\label{11jun14e5}
|D^{\alpha}(v_i^{k}(y)-V_i^k(\epsilon_ky)-2\log \epsilon_k-\Phi_i^k(y))|\le C\epsilon_k^2(1+|y|)^{4-m-|\alpha |+\delta}.
\end{equation}
where
$$\Phi_i^k(y)=\epsilon_k(G_{1,i}^k(r)\cos\theta+G_{2,i}^k(r)\sin \theta) $$
with
\begin{equation}\label{11jun14e6}
|G_{t,i}^k(r)|\le Cr(1+r)^{2-m+\delta}\quad t=1,2.
\end{equation}
\end{thm}

Note that $\Phi^k=(\Phi_1^k,..,\Phi_n^k)$ denotes the projection of $v_i^k$ onto $span\{\sin \theta, \cos\theta\}$. i.e.
\begin{equation}\label{11aug11e1}
\Phi^k_i(r\cos\theta,r\sin\theta)=\epsilon_k(G_{1,i}^k(r)\cos \theta+G_{2,i}^k(r)\sin\theta), \quad i\in I
\end{equation}
with $G_{t,i}^k(r)$ ($t=1,2$) satisfying some ordinary differential equations to be specified later.
The estimate for $|\alpha |=1$ follows from standard gradient estimate for elliptic equations.

\medskip

 Theorem \ref{cruexp} does not distinguish $m<4$ or $m=4$. But using Theorem \ref{cruexp} in more careful computation for $m<4$ and $m=4$ gives rise to more accurate results as follows: Here it is important to observe that $m-2<2$ if $m<4$. If we use
 $$m_i^k=\frac 1{2\pi}\int_{B(p,\delta)}\rho_i^k h_i^k(x)e^{u_i^k}dV_g, \quad m_k=\min\{m_1^k,...,m_n^k\}. $$
 Clearly $m_k\to m\in (2,4)$.

\begin{thm}\label{expthm1}
Suppose $m<4$, then for
$|y|\le \frac{\delta_0}2\epsilon_k^{-1}$ and $i\in I$
\begin{eqnarray}\label{mle4e1}
&&|D^{\alpha}\big (v_i^k(y)-(V_i^k(\epsilon_ky)+2\log \epsilon_k)-\Phi^k_i(y)\big )|\\
&\le &
C\epsilon_k^2(1+|y|)^{4-m_k-|\alpha|}\log(2+|y|).
\quad |\alpha |=0,1 \nonumber
\end{eqnarray}
 Moreover
$G_{t,i}^k$ ($t=1,2, i\in I$) satisfy
\begin{equation}\label{11jun9e1}
|G_{t,i}^k(r)|\le Cr(1+r)^{2-m_k}\quad 0<r<\epsilon_k^{-1}.
\end{equation}
\end{thm}

Note that Theorem \ref{expthm1} is slightly stronger than Theorem 4.2 of \cite{lin-zhang-jfa} because the latter has a logarithmic term. The reason is in the context of Theorem \ref{expthm1}, the function $v_i^k$ agrees with its approximation at the origin. Theorem \ref{expthm1} can be proved just like Theorem 4.2 in \cite{lin-zhang-jfa}.

\begin{thm}\label{expthm2}
If $m=4$ and $|m_i^k-4|\le C/\log \epsilon_k^{-1}$ for all $i\in I$, then we have, for
$|y|\le \frac{\delta_0}2\epsilon_k^{-1}$ and $i\in I$
\begin{eqnarray}\label{m4e1}
&&| D^{\alpha}(v_i^k(y)-(V_i^k(\epsilon_ky)+2\log \epsilon_k)-\Phi_i^{k}(y)\big )|\\
&\le &
C\epsilon_k^2(1+|y|)^{-|\alpha |}(\log (2+|y|))^2.
\quad |\alpha |=0,1, \nonumber
\end{eqnarray}
where $\Phi^{k}$ is of the form stated in (\ref{11aug11e1})
with $G_{t,i}^k$ ($t=1,2$) satisfying
\begin{equation}\label{11jun9e3}
|G_{t,i}^k(r)|\le Cr(1+r)^{-2},\quad 0<r<\epsilon_k^{-1},\quad i\in I.
\end{equation}
\end{thm}

\section{Rough estimates about bubbling magnitudes}
In this section, we will prove Theorems \ref{heightthm}, \ref{localmass} and \ref{mut-com}.
First in this section for simplicity we assume there are only two blowup points $p$ and $q$. The nature of analysis does not change if we have more blowup points.

Now we use Green's representation to describe the neighborhood of $p_k$. The expression of $u_i^k$ is
\begin{align*}
u_i^k(x)&=\bar u_i^k+\int_M G(x,\eta)\sum_{j=1}^n a_{ij} \rho_j^k h_j^k e^{u_j^k}dV_g\\
&=\bar u_i^k+\left(\int_{B(p_k,\delta)}+\int_{B(q_k,\delta)}+\int_{M\setminus (B(p_k,\delta)\cup B(q_k,\delta))}\right)G(x,\eta)\sum_{j=1}^n a_{ij} \rho_j^k h_j^k e^{u_j^k}dV_g \\
&=\bar u_i^k+I+II+III.
\end{align*}
Here
$\bar u_i^k=\int_M u_i^kdV_g$ and we use
$E=O(\epsilon_k^{m-2+\delta})$ to denote a very rough order of the error. If we use Theorem \ref{cruexp} in the evaluation, after cancellation we have
\begin{align*}
I&=\int_{B_{\delta}}(-\frac 1{2\pi}\log |\eta |+\gamma(x,\eta))\sum_{j=1}^n a_{ij}\tilde h_j^ke^{\tilde u_j^k}d\eta \\
&=-m_i^k \log |x|+2\pi m_i^k \gamma(x,p_k)+E.
\end{align*}
where $m_i^k= \frac{1}{2\pi}\sum_{j=1}^na_{ij}\sigma_j^k$ and $\sigma_i^k= \int_{B(p,\delta)}h_i^ke^{u_i^k}dV_g$, we use $\bar \sigma_i^k$, $\bar m_i^k$ to denote the integrations in $B(q_k,\delta)$. Here we note that in the evaluation of the integrals, the $\Phi_i^k$ terms disappear because of the symmetry of the domain.
Similarly the integration around $q_k$ gives
$$II=2\pi \bar m_i^kG(x,q_k)+E, \quad III=E. $$
Thus in the neighborhood of $p_k$ we have
\begin{equation}\label{ui-close-p}
u_i^k(x)=\bar u_i^k-m_i^k\log |x|+2\pi m_i^k\gamma(x,p_k)+2\pi \bar m_i^kG(x,q_k)+E
\end{equation}
in, say $B(p_k, 2\delta)\setminus B(p_k, \delta/2)$. If we use the approximation theorems to evaluate $u_i^k$ at $p_t^k$, it is easy to obtain
\begin{equation}\label{bad-est}
\bar u_i^k=(1-\frac{m_i^k}{2})M_k+O(1).
\end{equation}

Before more advanced estimates we first establish an elementary one:

\noindent{\bf Proof of Theorem \ref{heightthm}}:
From (\ref{bad-est}) we have
\begin{equation}\label{til-m-2}
(1-\frac{m_i^k}2)M_k=(1-\frac{\bar m_i^k}2)\bar M_k+O(1).
\end{equation}

Let
$$\lambda_k=M_k/\bar M_k, \quad \delta_i^k=O(1)/\bar M_k, $$
we have
\begin{equation}\label{tmk-mk}
(\frac{m_i^k}2-1)\lambda_k+\delta_k=(\frac{\bar m_i^k}2-1)+o(\epsilon_k^{\delta}), \quad \mbox{for some } \delta>0.
\end{equation}
It is established in \cite{linzhang1,lin-zhang-jfa} that $m^k=(m_1^k,...,m_n^k)$ satisfies
\begin{equation}\label{pi-mk}
\sum_{i=1}^n\sum_{j=1}^na^{ij}(\frac{m_i^k-2}2)(\frac{m_j^k-2}2)=\sum_{i=1}^n\sum_{j=1}^na^{ij}+E.
\end{equation}
and (\ref{pi-mk}) also holds for $\bar m^k=(\bar m_1^k,...,\bar m_n^k)$.
Thus using (\ref{tmk-mk}) in (\ref{pi-mk}) for $\bar m_k$, we have
$$ \sum_{i=1}^n\sum_{j=1}^na^{ij}(\frac{m_i^k-2}2\lambda_k+\delta_i^k)(\frac{m_j^k-2}2\lambda_k+\delta_j^k)=\sum_{i=1}^n\sum_{j=1}^na^{ij}+E, $$
which can be written as a quadratic expression of $\lambda_k$:
\begin{align*}
 &\lambda_k^2 \sum_{i=1}^n\sum_{j=1}^na^{ij}(\frac{m_i^k-2}2)(\frac{m_j^k-2}2)+2\lambda_k \sum_{i=1}^n\sum_{j=1}^na^{ij}(\frac{m_i^k-2}2)\delta_j^k+ \sum_{i=1}^n\sum_{j=1}^na^{ij}\delta_i^k\delta_j^k\\
 =& \sum_{i=1}^n\sum_{j=1}^na^{ij}+E/\bar M_k.
 \end{align*}
Let
$$B_k=\frac{ \sum_{i=1}^n\sum_{j=1}^na^{ij}(m_i^k-2)\delta_j^k}{ \sum_{i=1}^n\sum_{j=1}^na^{ij}}, \quad C_k=\frac{ \sum_{i=1}^n\sum_{j=1}^na^{ij}\delta_i^k\delta_j^k}{ \sum_{i=1}^n\sum_{j=1}^na^{ij}}. $$
Here we note that $ \sum_{i=1}^n\sum_{j=1}^na^{ij}>0$ because $(H2)$ requires $A$ to be invertible and $\sum_{j=1}^na^{ij}\ge 0$ for all $i$.
Then
$$\lambda_k^2+B_k\lambda_k+C_k=1+E/\bar M_k. $$
First it is obvious to observe that $\lim_{k\to \infty}\lambda_k=1$. Thus by $C_k=O(\bar M_k^{-2})$ and $B_k=O(\bar M_k^{-1})$
$$\lambda_k=\sqrt{1-C_k+\frac{B_k^2}4}-\frac{B_k}2=1-\frac{B_k}2-\frac 12(C_k-\frac{B_k^2}4)+O(\bar M_k^{-4}). $$
This verifies that $M_k-\bar M_k=O(1)$ which justifies $O(\epsilon_k^{m_k-2})=O(\bar \epsilon_k^{\bar m_k-2})$.
Theorem \ref{heightthm} is established. $\Box$

\medskip

By Theorem \ref{heightthm}, all the error terms above can be improved to $E=O(\epsilon_k^{m_k-2})$. Note that it is not $O(\epsilon_k^{m-2})$ yet, because the closeness of $m_k$ and $m$ is not derived yet. Another consequence of Theorem \ref{heightthm} is that
\begin{equation}\label{rho-close}
\frac{\sum_{j=1}^na_{ij}\rho_i^k}{2\pi}-m_i^k-\bar m_i^k=O(\epsilon_k^{m_k-2}).
\end{equation}
Indeed, integrating the equation for $u_i^k$ (which is (\ref{mainsys})), we have
$$\sum_{j=1}^n \int_M a_{ij}\rho_j^kh_j^ke^{u_j^k}=\sum_{j=1}^n a_{ij}\rho_j^k. $$
The integration of the left in $B(p,\delta)$ and $B(q,\delta)$ gives
$$m_i^k+\bar m_i^k+O(\epsilon_k^{m_k-2})=\sum_{j=1}^n a_{ij}\rho_j^k/(2\pi). $$
Thus (\ref{rho-close}) is verified.

Here we recall a theorem in \cite{lin-zhang-jfa} about location of blowup points:

Let $p_t^k$ be blowup points described as before. Then at each blowup point $p_t^k$, let $\phi_{it}^k$ be the harmonic function that eliminates the oscillation of
$\tilde u_i^k$ on $\partial B(\tilde p_t^k,\delta)$ for $\delta>0$ small. Then it is proved in \cite{lin-zhang-jfa} that

\begin{thm}\label{location}
If $m<4$
\begin{equation}\label{10050404}
|\sum_{i=1}^n\bigg (\partial_l (\log h_i^k)(p_t^k)+\partial_l \phi_{it}^k(p_t^k)\bigg )\rho_{it}^k|
\le C\epsilon_k^{m_k-2},\quad l=1,2,
\end{equation}
where $C$ is independent of $k$. On the other hand, if $m=4$
\begin{equation}\label{11jun28e1}
|\sum_{i=1}^n\bigg (\partial_l (\log h_i^k)(p_t^k)+\partial_l \phi_i^k(p_t^k)\bigg )\rho_{it}^k|
\le C\epsilon_k^2\log \epsilon_k^{-1},\quad l=1,2,
\end{equation}
where $\rho_{it}^k=\int_{B(p_t^k,\delta)}\rho_i^kh_i^ke^{u_i^k}dV_g$.
\end{thm}

In $B(p_k,\delta)$, by the definition of $\tilde u_i^k$ in (\ref{util-k}) and the estimate of $u_i^k$ in (\ref{ui-close-p}) we now have
\begin{align}\label{exp-p}
\tilde u_i^k(x)=\bar u_i^k-m_i^k\log |x|+2\pi m_i^k\gamma(x,p_k)+2\pi \bar m_i^k G(x,q_k)\\
+\log \rho_i^k+\log h_i^k(p_k)-f_i^k(x)+O(\epsilon_k^{m_k-2}).\nonumber
\end{align}

In this neighborhood, $\tilde u_i^k$ is of the form
$$\tilde u_i^k(x)=V_i^k(x)+\phi_i^k(x)+O(\epsilon_k^{m_k-2}),\quad x\in B(0,\delta/2)\setminus B(0,\delta/8) $$
where $\phi_i^k$ be the harmonic function on $B(p_k,\delta)$ that eliminates the oscillation of $\tilde u_i^k$:
$$\Delta \phi_i^k(x)=0,\quad \mbox{in}\quad B(p_k,\delta),\quad \phi_i^k(x)=\tilde u_i^k(x)-\frac{1}{2\pi\delta}\int_{\partial B(p_k,\delta)}\tilde u_i^k. $$
Note that we have used the fact that the first order terms $\Phi_i^k(x)=O(\epsilon_k^{m_k-2})$ when $x\sim 1$.
 Notice that by (\ref{rho-close}) and the fact that $q_k$ is not in $B(p_k,\delta)$,
\begin{align*}
&\Delta_g(2\pi m_i^k\gamma(x,p_k)+2\pi \bar m_i^kG(x,q_k)-f_i^k)\\
=&2\pi m_i^k+2\pi \bar m_i^k-\sum_{j=1}^n a_{ij}\rho_j^k=O(\epsilon_k^{m_k-2})
\end{align*}
Thus
\begin{equation}\label{who-har}
\Delta(2\pi m_i^k\gamma(x,p_k)+2\pi \bar m_i^kG(x,q_k)-f_i^k)=O(\epsilon_k^{m_k-2}).
\end{equation}
By the definition of $\tilde u_i^k$ in (\ref{util-k}) and the comparison of (\ref{ui-close-p}) and (\ref{who-har}) we have
\begin{align}\label{phi-2-b}
\phi_i^k(x)=2\pi m_i^k(\gamma(x,p_k)-\gamma(p_k,p_k))+2\pi \bar m_i^k(G(x,q_k)-G(p_k,q_k))\\
-f_i^k+O(\epsilon_k^{m_k-2}). \nonumber
\end{align}
 So if we rewrite the expression of $\tilde u_i^k$ as
\begin{align*}
\tilde u_i^k(x)&=\bar u_i^k-m_i^k\log |x|+\phi_i^k(x)\\
&+2\pi m_i^k \gamma(p_k,p_k)+2\pi \bar m_i^k(G(p_k,q_k))+\log (\rho_i^k h_i^k(p_k))+O(\epsilon_k^{m_k-2}).
\end{align*}
then we see that for $x\in B(0,\delta)\setminus B(0,\delta/8)$,
\begin{align}\label{vik-11}
&V_i^k(x)
=-m_i^k\log |x|+\bar u_i^k\\
&+2\pi m_i^k\gamma(p_k, p_k)+2\pi \bar m_i^k G(p_k,q_k)+\log (\rho_i^k h_i^k(p_k))+O(\epsilon_k^{m_k-2}). \nonumber
\end{align}

Similarly around $q_k$ we  have
\begin{align}\label{t-vik}
&\tilde V_i^k(x)
=-\bar m_i^k \log |x|+\bar u_i^k\\
&+2\pi \bar m_i^k \gamma(q_k,q_k)+2\pi m_i^k G(q_k, p_k)+\log (\rho_i^k h_i^k(q_k))+O(\epsilon_k^{m_k-2}). \nonumber
\end{align}

Let
$$M_k=\max_{i\in I}\max_{x}\tilde u_i^k(x),\quad \mbox{in}\quad B(p_k,\delta),$$
and $\bar M_k$ be the corresponding maximum in $B(q_k,\delta)$.
Then it is proved in \cite{linzhang2} that
$$M_k-\tilde u_i^k(p_k)=O(1).$$
We shall use the following notations:
$$D_i^k=\frac 1{2\pi}\int_{\mathbb R^2}\sum_{j=1}^n a_{ij}e^{V_j^k},\quad \alpha_i^k=M_k-\tilde u_i^k(p_k). $$
$\bar D_i^k$, $\bar \alpha_i^k$ can be understood similarly.

In order to obtain accurate estimate of $|m_i^k-\bar m_i^k|$, we first derive a simple fact about global solutions of Liouville system.

\begin{lem}\label{global-exp}
Let $U=(U_1,...,U_n)$ be global solution of
$$\Delta U_i+\sum_{j=1}^n a_{ij}e^{U_j}=0,\quad \mbox{in}\quad \mathbb R^2, \quad \int_{\mathbb R^2}e^{U_i}<\infty, \quad U_i \,\,\, \mbox{is radial} $$
and suppose $ \max_{i\in I}U_i(0)=0$.
Then
\begin{equation}\label{U-global}
U_i(r)=-m_i\log r+D_i-\alpha_i-\sum_{j=1}^n\frac{a_{ij}}{(m_j-2)^2}e^{D_j-\alpha_j}r^{2-m_j}+O(r^{2-m-\delta}).
\end{equation}
where $m_i=\frac{1}{2\pi}\int_{\mathbb{R}^2}\sum_{j=1}^n a_{ij}e^{U_j(x)}dx$,  $\alpha_i=-U_i(0)$
 and $D_i= \int_0^{\infty}\log r\sum_{j=1}^n a_{ij}e^{U_j(r)}rdr.$.
\end{lem}

\noindent{\bf Proof of Lemma \ref{global-exp}:}

It is easy to see that
$$U_i(x)=-\frac 1{2\pi}\int_{\mathbb R^2}\log |x-y|\sum_{j=1}^n a_{ij}e^{U_j}dy+c_i. $$
Recall that $U_i(0)=-\alpha_i$, then
$$-\alpha_i=-\int_0^{\infty}\log r\sum_{j=1}^n a_{ij}e^{U_j(r)}rdr+c_i. $$
and
\begin{align}\label{U-exp}
U_i(x)&=-\frac{1}{2\pi}\int_{\mathbb R^2}(\log |x-y|-\log |x|)\sum_{j=1}^n a_{ij}e^{U_j(y)}dy+D_i-\alpha_i-m_i\log |x|\\
&=-m_i\log |x|+D_i-\alpha_i+O(|x|^{-\delta}).\nonumber
\end{align}
for some $\delta>0$.
This expression gives
$$e^{U_i(r)}=r^{-m_i}e^{D_i-\alpha_i}+O(r^{-m-\delta}). $$
Then we use this in the ode that $U_i$ satisfies:
$$U_i''(r)+\frac 1rU_i'(r)=-\sum_{j=1}^n a_{ij}e^{U_j(r)},\quad 0<r<\infty. $$
Here we use the fact that
$$\lim_{r\to \infty} rU_i'(r)=-m_i. $$
Thus
\begin{align*}
-m_i-rU_i'(r)=-\sum_{j=1}^na_{ij}\int_r^{\infty}se^{U_j(s)}ds\\
=-\sum_{j=1}^n a_{ij}\frac{e^{D_j-\alpha_j}}{m_j-2}r^{2-m_j}+O(r^{2-m-\delta}).
\end{align*}
Then we have
$$U_i'(r)=-\frac{m_i}r+\sum_{j=1}^n\frac{a_{ij}}{m_j-2}e^{D_j-\alpha_j}r^{1-m_j}+O(r^{1-m-\delta}). $$
After integration and using (\ref{U-exp}) we see that (\ref{U-global}) holds.
Lemma \ref{global-exp} is established. $\Box $

\medskip

The main result in this section is:

\begin{prop}\label{really-close-m} If $m<4$,
$$|m_i^k-\bar m_i^k|\le C\delta_0^{4-m}\epsilon_k^{m_k-2},\quad i\in I. $$
\end{prop}

\noindent{\bf Proof of Proposition \ref{really-close-m}:}
Let $V_i^k$ be the sequence of global solutions that approximate $\tilde u_i^k$ around $p_k$, and in a neighborhood centered at $p_k$,
 $V_i^k(0)$ agrees with $\tilde u_i^k$ at $p_k$. $\bar V_i^k$ is understood as the first approximation around $q_k$ and we have
We use $m_{iv}^k$ to denote
$$m_{iv}^k=\frac 1{2\pi}\int_{\mathbb R^2}\sum_{j=1}^n a_{ij}e^{V_j^k}dx, \quad i=1,..,n, $$
and $\bar m_{iv}^k$ is for $\bar V_i^k$. Note that a rough estimate of $m_i^k$ based on Theorem \ref{expthm1} gives
$$m_i^k-m_{iv}^k=O(\epsilon_k^{m_k-2}). $$

If we use
$$U_i(y)=V_i^k(x)+2\log \epsilon_k,\quad \epsilon_k=e^{-\frac 12M_k},\quad x=\epsilon_k y. $$
Then the expansion of $V_i^k$ for $|x|$ is
\begin{align}\label{model-u}
V_i^k(x)&=-m_{iv}^k\log |x|-\frac{m_{iv}^k-2}{2}M_k+D_i-\alpha_i\\
&-\sum_{j=1}^n\frac{a_{ij}}{(m_{jv}^k-2)^2}e^{D_j-\alpha_j}\epsilon_k^{m_{jv}^k-2}|x|^{2-m_{jv}^k}+O(\epsilon_k^{m-2+\delta}). \nonumber
\end{align}
$V^k=(V_i^k,....,V_n^k)$ is the sequence of global functions that serves as the first term in the expansion of $\tilde u^k$ around $p_k$. Similarly
$$\bar V_i^k(x)=-\bar m_{iv}^k\log |x|-\frac{\bar m_{iv}^k-2}{2}\bar M_k+\bar D_i-\bar \alpha_i+O(\epsilon_k^{m_k-2})|x|^{2-\bar m_{iv}^k}$$

Since both $V^k=(V_1^k,...,V_n^k)$ and $\bar V^k=(\bar V_1^k,...,\bar V_n^k)$ are radial and satisfy the same Liouville system. The dependence on initial condition gives
\begin{equation}\label{key-1}
|V_i^k(x)-(\bar V_i^k(\eta x)+2\log \eta)|\le C\sum_{i=1}^n |\alpha_i-\bar \alpha_i|  \quad\mbox{in}\quad B(0,R).
\end{equation}
where $R>0$ is a constant, $\eta$ is chosen to make the heights equal, in this case $2\log \eta=M_k-\bar M_k$. We also note that one of $\alpha_i=\bar \alpha_i$. Here we invoke an important result in \cite{linzhang1}. Suppose $\alpha_1=0$, the mapping from $(\alpha_2,...,\alpha_{n})$ to $(\sigma_2,...,\sigma_{n})$ is invertible. In fact the following matrix
$$\mathbb M=\left(\begin{array}{ccc}
\partial_{\alpha_2}\sigma_2 & ... & \partial_{\alpha_n}\sigma_2 \\
\vdots & \vdots & \vdots \\
\partial_{\alpha_2}\sigma_n & ... & \partial_{\alpha_n}\sigma_n
\end{array}
\right) $$
is invertible. Note that in this proposition we use
$$\sigma_i^k=\sum_{j=1}^na^{ij}m_{jv}^k, \quad \bar \sigma_i^k=\sum_{j=1}^na^{ij}\bar m_{jv}^k. $$

 Thus (\ref{key-1}) can be written as
\begin{equation}\label{key-2}
|V_i^k(x)-(\bar V_i^k(\eta x)+2\log \eta)|\le C\sum_{i=2}^n|\sigma_i^k-\bar \sigma_i^k|.
\end{equation}
By the expansion of $\bar V_i^k$, we find that
\begin{equation}
 \bar V_i^k(\eta x)+2\log \eta=-\frac{\bar m_{iv}^k-2}2M_k -\bar m_{iv}^k\log |x|+\bar D_i-\bar \alpha_i+O(\epsilon_k^{m_k-2}).
\end{equation}
The difference between $V_i^k$ and $\bar V_i^k(\eta x)+2\log \eta$ gives
\begin{align*}
&V_i^k(x)-(\bar V_i^k(\eta x)+2\log \eta)\\
=&(\bar m_{iv}^k-m_{iv}^k)\log |x|+\frac{\bar m_{iv}^k-m_{iv}^k}2 M_k +D_i-\bar D_i+\bar \alpha_i-\alpha_i+O(\epsilon_k^{m_k-2}).
\end{align*}
By the dependence of initial condition and fast decay of $V^k$ and $\bar V^k$, we have
$$|D_i-\bar D_i|\le C\sum_{i=1}^n|\alpha_i-\bar \alpha_i|\le C\sum_{i=2}^n|\sigma_i^k-\bar \sigma_i^k|. $$
And by the one to one correspondence between $(\alpha_1,...,\alpha_n)$ to $(\sigma_2,...,\sigma_n)$ and the smoothness of the mapping we also obtain
$$\sum_{i=1}^n|\alpha_i-\bar \alpha_i|\le C\sum_{i=2}^n|\sigma_i^k-\bar \sigma_i^k|. $$
Combining these estimates we have
\begin{equation}\label{imp-1}
 |\frac{\bar m_{iv}^k-m_{iv}^k}2 M_k |\\
\le  C\sum_{j=2}^n|\sigma_{j}^k-\bar \sigma_{j}^k|+O(\epsilon_k^{m_k-2}) \quad i=1,...,n.
\end{equation}
  After that we multiply $a^{ij}$ with summation on $i$ and take summation on $j$, then we have
\begin{equation}\label{imp-2}
\sum_{j=1}^n|\bar \sigma_j^k-\sigma_j^k |\le \frac C{M_k}\sum_{l=2}^n|\sigma_l^k-\bar \sigma_l^k|+O(\epsilon_k^{m_k-2}).
\end{equation}
We thus obtain the following important closeness result:
\begin{equation}\label{imp-4}
\sigma_i^k-\bar \sigma_i^k=O(\epsilon_k^{m_k-2})/M_k,\quad i=1,...,n.
\end{equation}
Thus we have proved that $m_{iv}^k-\bar m_{iv}^k=O(\epsilon_k^{m_k-2}/M_k)$, using the expansion of $\tilde u_i^k$ in the calculation of
$\int_{B(p_k,\delta_0)}\rho_i^kh_i^ke^{u_i^k}dV_g$. It is easy to see that
$$|m_i^k-m_{iv}^k|\le C\delta_0^{m_k-4}\epsilon_k^{m_k-2}, \quad |\bar m_i^k-\bar m_{iv}^k|\le C\delta_0^{4-m}\epsilon_k^{m_k-2}. $$
Thus
$$|m_i^k-\bar m_i^k|\le C\delta_0^{4-m}\epsilon_k^{m_k-2}. $$
Proposition \ref{really-close-m} is established. $\Box$

\medskip

\begin{rem} By the same argument for the case $m=4$, we also have
\begin{equation}\label{4-g-close}
m_{iv}^k-\bar m_{iv}^k=O(\epsilon_k^{m_k-2}/M_k).
\end{equation}
 Here $m_k\to 4$. Later we shall see $m_k$ can be replaced by $4$. Correspondingly
\begin{equation}\label{4-close}
m_i^k-\bar m_i^k=O(\epsilon_k^2\log \frac{1}{\epsilon_k}),
\end{equation}
but the leading term (of the order $O(\epsilon_k^2\log \frac{1}{\epsilon_k})$ ) can be identified as a local term that involves curvature at the blowup point.
\end{rem}

In the more general situation of $N$ bubbles, if we use $m_{it}^k$ to denote the energy in $B(p_t^k,\delta_0)$, $m_{itv}^k$ to denote the energy of the global sequence as the first term in the approximation,  we have, for $m<4$:
\begin{equation}\label{m-close-m}
m_{it}^k-m_{is}^k=O(\epsilon_k^{m_k-2}), \quad s\neq t, \quad m<4
\end{equation}
\begin{equation}\label{m-close-1}
m_{itv}^k-m_{isv}^k=O(\epsilon_k^{m_k-2}/\log \frac{1}{\epsilon_k}), \quad s\neq t, \quad m<4.
\end{equation}
For $m=4$, (\ref{m-close-1}) also holds, and $m_{it}^k-m_{is}^k=O(\epsilon_k^{m_k-2}\log \epsilon_k^{-1})$. The following lemma gives an estimate of $\rho_i^k-\rho$, which determines as a consequence that
$$O(\epsilon_k^{m_k-2})=O(\epsilon_k^{m-2}) \quad \mbox{if}\quad m\in (2,4)$$
$$O(\epsilon_k^{m_k-2}\log \frac{1}{\epsilon_k})=O(\epsilon_k^2\log \frac{1}{\epsilon_k}),\quad \mbox{if}\quad m=4. $$

\begin{lem}\label{rho-close-1} If $m<4$
$$\rho_i^k-\rho_i=O(\epsilon_k^{m-2}),\quad i\in I.$$
\end{lem}

\noindent{\bf Proof of Lemma \ref{rho-close-1}:} Recall that $\rho\in \Gamma_N$.
Let
$$\rho_{it}^k=\int_{B(p_t,\delta)}\rho_ih_i^ke^{u_i^k}dV_g,\quad t=1,...,N,\quad E_i^k=\rho_i^k-\sum_{t=1}^N\rho_{it}^k. $$
Here $\delta>0$ is small so that bubbling disks are disjoint. Clearly from $\int_{M}h_i^ke^{u_i^k}=1$ and (\ref{ui-bar}) we have
$$E_i^k=O(\epsilon_k^{m_k-2}),\quad \rho_i^k=\sum_{i=1}^N\rho_{it}^k+E_i^k. $$
Let $\sigma_i=\frac{\rho_i}{2\pi N}$ and $\sigma_{it}^k=\frac{\rho_{it}^k}{2\pi}$. Then we write $\sigma_i=\sigma_{it}+s_{it}^k$. We have known that $s_{it}^k=o(1)$ as $k\to \infty$. Since $\rho\in \Gamma_N$ we know
$$ \sum_{i=1}^n\sum_{j=1}^na_{ij}\sigma_i\sigma_j-4\sum_{i=1}^n\sigma_i=0. $$
On the other hand the Pohozaev identity around $p_t$ gives
$$ \sum_{i=1}^n\sum_{j=1}^na_{ij}\sigma_{it}^k\sigma_{jt}^k-4\sum_{i=1}^n \sigma_{it}^k=O(\epsilon_k^{m_k-2}). $$
The difference between these two equations gives
\begin{equation}\label{s-close}
\sum_{i=1}^n 2(m_i-2)s_{it}^k+ \sum_{i=1}^n\sum_{j=1}^na_{ij}s_{it}^ks_{jt}^k=O(\epsilon_k^{m_k-2}),\quad t=1,...,N.
\end{equation}
Taking the sum of $N$ equations we have
$$\sum_{i=1}^n2(m_i-2)(\sum_{t=1}^Ns_{it}^k)+\sum_{t=1}^Na_{ij}s_{it}^ks_{jt}^k=O(\epsilon_k^{m_k-2}). $$
From Proposition \ref{really-close-m} we know that the difference between any two $s_{it}^k$ is $O(\epsilon_k^{m_k-2})$, $\sum_{t=1}^Ns_{it}^N=Ns_{i1}^k+O(\epsilon_k^{m_k-2})$.
Thus we have
$$N\sum_{i=1}^ns_{i1}^k+\sum_{t=1}^N \sum_{i=1}^n\sum_{j=1}^na_{ij}s_{it}^ks_{jt}^k=O(\epsilon_k^{m_k-2}).$$
By the assumption
$\frac{\rho_i^k-\rho_i}{\rho_j^k-\rho_j}\sim 1 $ we have $s_{i1}^k/s_{j1}^k\sim 1$ for all $i,j$ because $\sum_{t=1}^N s_{it}^k=\frac 1{2\pi}(\rho_i^k-\rho_i)$.
Thus
$$s_{it}^k=O(\epsilon_k^{m_k-2}),\quad t=1,...,N.$$
 Thus Lemma \ref{rho-close-1} is established. $\Box$\par
Now we are ready to prove Theorem \ref{localmass}.\par

\noindent{\bf Proof of  Theorem \ref{localmass}:}
Due to the closeness between $m_i^k$ and $m$ (which is of the order $O(\epsilon_k^{m_k-2})$ we shall use $O(\epsilon_k^{m-2})$ instead of $O(\epsilon_k^{m_k-2})$ in Proposition \ref{really-close-m}.
  By a similar argument, based on (\ref{4-close}) we clearly have
$$\mbox{ For } m=4, \quad \rho_i^k-\rho_i=O(\epsilon_k^{2}\log \frac{1}{\epsilon_k}), $$
which gives $O(\epsilon_k^{m_k-2})=O(\epsilon_k^2)$.   We thus have
$$m_{i,t}^k-m_{i,s}^k=O(\varepsilon_k^2\log\frac{1}{\varepsilon_k}),\,\, t\not=s,\,\,  s,t\in\{1,\cdots, N\},\,\,  i\in I.$$
Recall that $\sigma_{i,t}^k=2\pi\sum_{j\in I}a^{ij}m_{i,t}^k$ and $A=(a_{ij})$ is invertible. This theorem is proved.
   $\Box$

Because of (\ref{m-close-m}) the evaluation of $u_i^k(x)$ away from bubbling disks now becomes:
$$u_i^k(x)=\bar u_i^k+2\pi m_i^k\sum_{l=1}^N G^*(x,p_l^k)+O(\epsilon_k^{m-2}), \quad x\in M\setminus (\cup B(p_t^k,\delta_0)), \quad m<4,$$
$$u_i^k(x)=\bar u_i^k+2\pi m_i^k\sum_{l=1}^N G^*(x,p_l^k)+O(\epsilon_k^{2}\log \frac{1}{\epsilon_k}), \quad x\in M\setminus (\cup B(p_t^k,\delta_0)), \quad m=4,$$
for $\delta_0>0$ small.  We recall the notation
$$G^*(p_t^k,p_l^k)=\left\{\begin{array}{ll} \gamma(p_t^k,p_t^k),\quad \mbox{if}\quad t=l, \\
G(p_t^k,p_l^k), \quad \mbox{if}\quad t\neq l.
\end{array}
\right.
$$
With the information available we are in the position to prove Theorem \ref{mut-com}.

\medskip

\noindent{\bf Proof of Theorem \ref{mut-com}:}  In $B(p_t^k,\delta_0)$ we use $V_t^k=(V_{1t}^k,...,V_{nt}^k)$ to denote the sequence of global solutions as the first term in the approximation. In the context of multiple bubbles, $V_{it}^k$ has two expressions:
$$V_{it}^k=-m_{it}^k\log |x|+\bar u_i^k+2\pi m_i^k\sum_{l=1}^N G^*(p_t^k,p_l^k)+\log (\rho_i^kh_i^k(p_t^k))+O(\epsilon_k^{m-2}), $$
and
$$V_{it}^k=-m_{itv}^k\log |x|-\frac{m_{itv}^k-2}2 M_{k,t}+D_{it}^k-\alpha_{it}^k+O(\epsilon_k^{m-2}), $$
where $m_{itv}^k=\sum_{j=1}^na_{ij}\sigma_{itv}^k$, $\sigma_{itv}^k=\frac 1{2\pi}\int_{\mathbb R^2}e^{V_{it}^k}$.

By comparing the two expressions of $V_{it}^k$ we have

\begin{align*}-m_{it}^k\log |x|+\bar u_i^k+2\pi m_{it}^k\sum_{l=1}^N G^*(p_t^k,p_l^k)+\log (h_i^k(p_t^k)\rho_i^k)\\
=-m_{itv}^k\log |x|-\frac{m_{itv}^k-2}2M_{k,t}+D_{it}^k-\alpha_{it}^k+O(\epsilon_k^{m-2}).
\end{align*}
where all the notations are clearly understandable under this context.
Thus, for $t\neq s$, based on the two different expression of $\bar u_i^k$, we have
$$\sum_{l=1}^N 2\pi m_i^k G^*(p_t^k,p_l^k)+\log h_i^k(p_t^k)=\sum_{l=1}^N 2\pi m_i^k G^*(p_s^k,p_l^k)+\log h_i^k(p_s^k)+O(\epsilon_k^{m-2}), $$
where we have used $m_{ivt}^k-m_{ivs}^k=O(\epsilon_k^{m-2}/\log \frac{1}{\epsilon_k})$ and used $m_i^k$ to replace all $m_{it}^k$. (\ref{mut-re}) is verified.
(\ref{mut-re-2}) can be verified similarly. Theorem \ref{mut-com} is established. $\Box$

\medskip

\medskip

The following expression of $\bar u_i^k$ will be used
\begin{align}\label{ui-bar}
\bar u_i^k&=(1-\frac{m_i^k}2)M_{k,t}-\log (\rho_i^kh_i^k(p_t^k))-2\pi m_i^k\sum_{l=1}^N G^*(p_t^k,p_l^k)\\
&+D_i^k-\alpha_i^k+O(\epsilon_k^{m_k-2}),\quad t=1,...,N, \nonumber
\end{align}
where we used $\alpha_i^k$ and $D_i^k$ to replace any $D_{it}^k$ and $\alpha_{it}^k$ for obvious reasons. One of $M_{k,t}$ is $M_k$ and the difference between any two
of these is bounded. In fact, for $t\neq s$, the difference on equations (\ref{ui-bar}) for $t$ and $s$ gives
$$2\pi m_i(\sum_{l=1}^N (G^*(p_t^k,p_l^k)-G^*(p_s^k,p_l^k))+\log \frac{h_i^k(p_t^k)}{h_i^k(p_s^k)}=-\frac{m_i-2}2(M_{k,t}-M_{k,s})+O(\epsilon_k^{m-2}). $$
Equivalent form is
\begin{equation}\label{s-t-comp}
exp(2\pi m_i\sum_{l=1}^N (G^*(p_t^k,p_l^k)-G^*(p_s^k,p_l^k)))\, \frac{h_i^k(p_t^k)}{h_i^k(p_s^k)}=\frac{\epsilon_{k,t}^{m_i-2}}{\epsilon_{k,s}^{m_i-2}}+
O(\epsilon_k^{m-2}),
\end{equation}
where $\epsilon_{k,t}=e^{-\frac 12 M_{k,t}}$. Also, (\ref{ui-bar}) gives
\begin{equation}\label{e-uib}
e^{\bar u_i^k}=\epsilon_{k,t}^{m_i-2}\frac{e^{D_i^k-\alpha_i^k}}{\rho_i^kh_i^k(p_t^k)}e^{-2\pi m_i\sum_{l=1}^N G^*(p_t^k,p_l^k)}+O(\epsilon_k^{m-2+\delta}).
\end{equation}
\section{Proof of leading terms}

Now we are in the position to prove Theorem \ref{gammaNpneq}.

\medskip

\noindent{\bf Proof of Theorem \ref{gammaNpneq}:}

 Since $\int_M h_i^ke^{u_i^k}dV_g=1$ we write
$$\rho_i^k=\sum_{t=1}^N\int_{B(p_t^k,\delta_0)}\rho_i^kh_i^ke^{u_i^k}dV_g+\int_{M\setminus \cup_tB(p_t^k,\delta_0)}\rho_i^kh_i^ke^{u_i^k}dV_g
=\sum_{t=1}^N\rho_{it}^k+\rho_{ib}^k $$
where in local coordinates
$$
\rho_{it}^k=\int_{B(0,\delta_0)}\tilde h_i^ke^{\tilde u_i^k}d\eta.
$$
Let
$$I_1=\{i\in I;\quad \lim_{k\to \infty}m_i^k=m.\quad \}. $$
Based on Proposition \ref{really-close-m} for each $i\in I_1$, $m_i^k-m=O(\epsilon_k^{m-2})$.
If we use $V_i^k$ to be the leading term in the approximation of $\tilde u_i^k$ and $U_i^k$ be the scaled version of $V_i^k$, by
 (\ref{model-u}) we have (since $m<4$)
\begin{align}\label{rhoia}
\frac{1}{2\pi}\rho_{it}^k&=\frac 1{2\pi}\int_{B(0,\delta_0\epsilon_k^{-1})}\tilde h_i^k(0)e^{V_i^k(y)}dy+o(\delta_0)\epsilon_k^{m-2},\\
&=\sigma_{iv}^k-\frac{e^{D_i-\alpha_i}}{m-2}\epsilon_{k,t}^{m-2}\delta_0^{2-m}+E_{\delta_0},\quad i\in I_1.  \nonumber
\end{align}
where we use $\sigma_{iv}^k=\frac 1{2\pi}\int_{\mathbb R^2}e^{V_{it}^k}$ and $E_{\delta_0}$ to denote $o(\delta_0)\epsilon_k^{m-2}$, $\epsilon_{k,t}=e^{-\frac{M_{k,t}}2}$. We did not use $t$ in $\sigma_{iv}^k$ because the difference between any two of them is $O(\epsilon_k^{m-2}/M_k)$.
Now for $i\not \in I_1$ we have
\begin{equation}\label{rhoia1}
\frac 1{2\pi}\rho_{it}^k=\sigma_{iv}^k+E_{\delta_0},\quad i\not \in I_1
\end{equation}
and
\begin{equation}\label{rhobc}
|\rho_{ib}^k|=E_{\delta_0},\quad i\not \in I_1.
\end{equation}

Combining (\ref{rhoia}), (\ref{rhoia1}), (\ref{rhobc})  we have
\begin{align*}
&\sum_{i=1}^n \frac{4}{2\pi}\rho_{it}^k- \sum_{i=1}^n\sum_{j=1}^na_{ij}\frac{\rho_{it}^k}{2\pi}\frac{\rho_{jt}^k}{2\pi}\\
=&4\sum_{i=1}^n(\sigma_{iv} -\frac{e^{D_i-\alpha_i}}{m_i-2}\delta^{2-m_i}\epsilon_{k,t}^{m_i-2})\\
&- \sum_{i=1}^n\sum_{j=1}^na_{ij}(\sigma_{iv}^k-\frac{e^{D_i-\alpha_i}}{m_i-2}\delta^{2-m_i}\epsilon_{k,t}^{m_i-2})(\sigma_{jv}^k-\frac{e^{D_j-\alpha_j}}{m_j-2}\delta^{2-m_j}\epsilon_{k,t}^{m_j-2})+E_{\delta_0}\\
=&-\frac{4}{m-2}\sum_{i\in I_1}e^{D_i-\alpha_i}\delta_0^{2-m}\epsilon_{k,t}^{m-2}+2\sum_{i\in I_1} \sum_{j=1}^na_{ij}\sigma_{jv}^k\frac{e^{D_j-\alpha_j}}{m_j-2}\delta_0^{2-m_j}\epsilon_{k,t}^{m_j-2}+E_{\delta_0}\\
=&2\delta_0^{2-m}\epsilon_{k,t}^{m-2}\sum_{i\in I_1}e^{D_i-\alpha_i}+E_{\delta_0}.
\end{align*}
Note that $m_i= \sum_{j=1}^na_{ij}\sigma_{jv}^k+O(\epsilon_k^{m-2}/M_k)$.  For $i\in I_1$, using (\ref{ui-bar}) we have
\begin{align*}
\rho_{ib}^k&=\int_{M\setminus (\cup_t B(p_t^k,\delta_0))}\rho_i^kh_i^ke^{u_i^k}dV_g\\
=&\rho_i^ke^{\bar u_i^k}\int_{M\setminus \cup_t B(p_t^k,\delta_0)}h_i^ke^{2\pi m\sum_{t=1}^N G(x,p_t^k)}
+E_{\delta_0},\quad i\in I_1.
\end{align*}
Now we define $N$ open sets $\Omega_{t,\delta_0}$ such that they are mutually disjoint, each of them contains a bubbling disk and their union is $M$:
$$B(p_t^k,\delta_0)\subset \Omega_{t,\delta_0}, \quad \cup_{t=1}^N \overline{\Omega_{t,\delta_0}}=M,
\quad \Omega_{t,\delta_0}\cap \Omega_{s,\delta_0}=\emptyset, \,\, \forall t\neq s. $$
In each $\Omega_{t,\delta_0}$ we use (\ref{e-uib})) to
 write $\rho_{ib}^k$ as (for $i\in I_1$)
\begin{align}\label{rhoib-2}
&\rho_{ib}^k=\rho_i^ke^{\bar u_i^k}\sum_{t=1}^N\int_{\hat{\Omega}_{t,\delta_0}}h_i^ke^{2\pi m\sum_{l=1}^NG(x,p_l^k)}\\
&=\sum_{t=1}^N\int_{\hat{\Omega}_{t,\delta_0}}\epsilon_{k,t}^{m-2}\frac{h_i^k(x)}{h_i^k(p_t^k)}e^{D_i-\alpha_i}e^{2\pi m\sum_{l=1}^N (G(x,p_l^k)-G^*(p_t^k,p_l^k))}dV_g
\nonumber\\
&+E_{\delta_0},\quad i\in I_1,\nonumber
\end{align}
where $\hat{\Omega}_{t,\delta_0}= \Omega _{t,\delta_0}\setminus B(p_t^k,\delta_0)$. Now we put estimates together to have
\begin{align*}
&4\sum_{i=1}^n\frac{\rho_i^k}{2\pi N}-\sum_{i=1}^n\sum_{j=1}^na_{ij}\frac{\rho_i^k}{2\pi N}\frac{\rho_j^k}{2\pi N}\\
=&4\sum_{i=1}^n(\sum_{t=1}^N \frac{\rho_{it}^k}{2\pi N}+\frac{\rho_{ib}}{2\pi N })- \sum_{i=1}^n\sum_{j=1}^na_{ij}(\sum_{t=1}^N \frac{\rho_{it}^k}{2\pi N} +\frac{\rho_{ib}^k}{2\pi N})
(\sum_{s=1}^N\frac{\rho_{js}^k}{2\pi N}+\frac{\rho_{jb}^k}{2\pi N})\\
=&\frac{2}{N}\sum_{t=1}^N\sum_{i\in I_1}\delta_0^{2-m}\epsilon_{k,t}^{m-2}e^{D_i-\alpha_i}- 2(m-2)\sum_{i\in I_1}\frac{\rho_{ib}^k}{2\pi N}+E_{\delta_0}.
\end{align*}

Using (\ref{rhoib-2}) in the expression above we have
\begin{align}\label{leading-t}
&\qquad 4\sum_{i=1}^n\frac{\rho_i^k}{2\pi N}-\sum_{i=1}^n\sum_{j=1}^na_{ij}\frac{\rho_i^k}{2\pi N}\frac{\rho_j^k}{2\pi N}\\
&=\frac{2}{N}\sum_{i\in I_1}e^{D_i-\alpha_i}\sum_{t=1}^N\epsilon_{k,t}^{m-2}\bigg (\delta_0^{2-m}-\frac{(m-2)}{2\pi}\int_{\hat{\Omega}_{t,\delta_0}}
(\frac{h_i^k(x)}{h_i^k(p_t^k)}e^{2\pi m\sum_{l=1}^N(G(x,p_l^k)-G^*(p_t^k,p_l^k))}dV_g\bigg )\nonumber\\
&+E_{\delta_0} \nonumber
\end{align}

Since one of $\epsilon_{k,t}$ is $\epsilon_k$, say $\epsilon_{1,t}=\epsilon_k$, based on (\ref{s-t-comp}) we use
\begin{equation}\label{c-t-def}
c_t=\epsilon_{k,t}^{m-2}/\epsilon_k^{m-2}=\frac{h_i^k(p_t^k)e^{2\pi m\sum_{l=1}^NG^*(p_t^k,p_l^k)}}{h_i^k(p_1^k)e^{2\pi m\sum_{l=1}^N G^*(p_1^k,p_l^k)}},
\end{equation}
Here we note that (\ref{s-t-comp}) implies that for $i\in I_1$, $c_t$ is independent of $i\in I_1$. Then we use $D+o(1)$ to represent
\begin{align}\label{notation-D}
&D+o(1):=\\
&\sum_{i\in I_1}e^{D_i-\alpha_i}\sum_{t=1}^Nc_t\bigg (\delta_0^{2-m}-\frac{(m-2)}{2\pi}\int_{\hat{\Omega}_{t,\delta_0}}
(\frac{h_i^k(x)}{h_i^k(p_t^k)}e^{2\pi m\sum_{l=1}^N(G(x,p_l^k)-G^*(p_t^k,p_l^k))}dV_g\bigg ). \nonumber
\end{align}
and the leading term of $4\sum_{i=1}^n\frac{\rho_i^k}{2\pi N}-\sum_{i=1}^n\sum_{j=1}^n\frac{\rho_i^k}{2\pi N}\frac{\rho_j^k}{2\pi N}$ is written as
$$4\sum_{i=1}^n\frac{\rho_i^k}{2\pi N}-\sum_{i=1}^n\sum_{j=1}^n\frac{\rho_i^k}{2\pi N}\frac{\rho_j^k}{2\pi N}=\frac 2N\epsilon_k^{m-2}(D+o(1)), $$
where $o(1)$ stands for the infinitesimal quantity when $\delta_0\to 0$. Here it is important to observe that $D$ is involved with integration on the whole manifold.
Theorem \ref{gammaNpneq} is established. $\Box$

\bigskip

\noindent{\bf Proof of Theorem \ref{main-3}:}

Around each $p_t^k$, an extension of (\ref{phi-2-b}) can be easily determined to be
\begin{align}\label{phi-pt}
\phi_i^k(x)&=2\pi m_i(\gamma(x,p_t^k)-\gamma(p_t^k,p_t^k))\\
&+\sum_{l\neq t} 2\pi m_i(G^*(x,p_l^k)-G^*(p_t^k,p_l^k))-f_i^k(x)+E_1 \nonumber
\end{align}
where $E_1=O(\epsilon_k^{m-2})$ if $m<4$ and is $O(\epsilon_k^2\log \frac{1}{\epsilon_k})$ if $m=4$.
Correspondingly,
\begin{equation}\label{phi-ptd}
\nabla\phi_i^k(p_t^k)=2\pi  m_i\sum_{l=1}^N\nabla_1G^*(p_t^k,p_l^k)+E_1.
\end{equation}
With these notations, (\ref{11july13e7}) and (\ref{11july13e8}) follow immediately. Theorem \ref{main-3} is established. $\Box$

\medskip

Finally we prove Theorem \ref{gammaNpeq}.

\medskip

\noindent{\bf Proof of Theorem \ref{gammaNpeq}:}

$$\rho_i^k=\sum_{t=1}^N\int_{B(p_t^k,\delta_0)}\rho_i^k h_i^k e^{u_i^k} dV_g +\int_{M\setminus \cup_tB(p_t^k,\delta_0)} \rho_i^k h_i^k e^{u_i^k} dV_g. $$
We continue to use the notation $\rho_{it}^k$ and $\rho_{ib}^k$.
By (\ref{bad-est}) and Theorem \ref{expthm2} the second integral is $O(\epsilon_k^2)$, this is the same as the computation for the single equation \cite{ChenLin1}.

Now we use the expansion of bubbles to compute each $\rho_{it}^k$. By the expansion of $\tilde u_i^k$ around $p_t^k$, we have
\begin{align}\label{11july5e1}
&\rho_{it}^k= \int_{B(p_t^k,\delta_0)}\rho_i^k h_i^k e^{u_i^k} dV_g
=\int_{B(0,\delta_0)}\tilde h_i^k e^{\phi_i^k} e^{\tilde u_i^k-\phi_i^k}d\eta\nonumber\\
&=\int_{B(0,\delta_0 \epsilon_{k,t}^{-1})} \rho_i^k h_i^k(p_t^k) e^{U_i^k(\eta)}d\eta+O(\epsilon_k^2)\\
&+\int_{B(0,\delta_0 \epsilon_{k,t}^{-1})}
\epsilon_{k,t}^2(\frac 14 \Delta \tilde h_i^k(0)+\frac 12\nabla \tilde h_i^k(0)\cdot \nabla \phi_i^k(0)+\frac 14|\nabla \phi_i^k(0)|^2)
|y|^2e^{U_i^k}d\eta.\nonumber
\end{align}
The first integral on the right hand side of the above is $O(\epsilon_k^2)$ different from the global solution in the approximation of $\tilde u_i^k$ around $p_t^k$.
So we use $\sigma_{ivt}^k$ to denote it. For $t\neq s$, from (\ref{4-g-close}) we see that
$$\sigma_{ivt}^k-\sigma_{ivs}^k=O(\epsilon_k^2/\log \frac{1}{\epsilon_k}). $$
 To evaluate the last term, we first use the
definition of the $\tilde h_i^k$ to have

$$\nabla\tilde h_i^k(0)\cdot \nabla \phi_i^k(0)=2\pi m_i\frac{\nabla h_i^k(p_t^k)}{h_i^k(p_t^k)}\cdot \sum_{l=1}^N\nabla_1G^*(p_t^k,p_l^k)+O(\epsilon_k^2\log \frac{1}{\epsilon_k}),$$
and
\begin{align}\label{h-lap}
\Delta \tilde h_i^k(0)=\frac{\Delta h_i^k(p_t^k)}{h_i^k(p_t^k)}-2K(p_t^k)+\sum_{j=1}^na_{ij}\rho_j+O(\epsilon_k^2\log \frac{1}{\epsilon_k})\nonumber\\
=\frac{\Delta h_i^k(p_t^k)}{h_i^k(p_t^k)}-2K(p_t^k)+8\pi N+O(\epsilon_k^2\log \frac{1}{\epsilon_k}).
\end{align}

Then we define $b_{it}^k$ as
\begin{align}\label{bit}
b_{it}^k&=e^{D_i-\alpha_i}\bigg (\frac 14\frac{\Delta h_i^k(p_t^k)}{h_i^k(p_t^k)}-K(p_t^k)+4\pi N\\
&+4\pi\frac{\nabla h_i^k(p_t^k)}{h_i^k(p_t^k)}\cdot \sum_{l=1}^N\nabla_1G^*(p_t^k,p_l^k)
+16\pi^2|\sum_{l=1}^N\nabla_1G^*(p_t^k,p_l^k)|^2\bigg ). \nonumber
\end{align}
With this $b_{it}^k$ we have
$$\frac{\rho_{it}^k}{2\pi}=\sigma_{iv}^k+b_{it}^k\epsilon_{k,t}^2\log \frac{1}{\epsilon_k}+O(\epsilon_k^2). $$
From (\ref{s-t-comp}) we see that $\epsilon_{k,t}$ can be replaced by $\epsilon_k$.
Consequently,
\begin{align}\label{11jun8e1}
&4\sum_{i=1}^n\frac{\rho_i^k}{2\pi N}- \sum_{i=1}^n\sum_{j=1}^na_{ij}\frac{\rho_i^k}{2\pi N}\frac{\rho_j^k}{2\pi N}\\
=&4\sum_{i=1}^n\sum_{t=1}^N\frac{\rho_{it}^k}{2N \pi}- \sum_{i=1}^n\sum_{j=1}^na_{ij}(\sum_{t=1}^N \frac{\rho_{it}}{2\pi N})(\sum_{s=1}^N \frac{\rho_{js}}{2\pi N})\nonumber \\
=&4\sum_{i=1}^n\sum_{t=1}^N(\frac{\sigma_{ivt}}N+\tilde \epsilon_k\frac{b_{it}^k}{N})-\sum_{i=1}^n\sum_{j=1}^n\sum_{t=1}^N\sum_{s=1}^Na_{ij}(\frac{\sigma_{ivt}}N+\tilde\epsilon_k\frac{b^k_{it}}N)(\frac{\sigma_{jv}}N+\tilde\epsilon_k\frac{b^k_{js}}N)\nonumber \\
=&4\sum_{i=1}^n\frac{\sigma_{iv}}N- \sum_{i=1}^n\sum_{j=1}^na_{ij}\frac{\sigma_{iv}\sigma_{jv}}N+4\sum_{i=1}^n\sum_{t=1}^N\tilde\epsilon_k\frac{b^k_{it}}N
-2\sum_{i=1}^n\sum_{j=1}^n\sum_{t=1}^Na_{ij}\sigma_{iv}\tilde\epsilon_k\frac{b^k_{jt}}N+O(\epsilon_k^2)
\nonumber \\
=&-4\epsilon_k^2\log \epsilon_k^{-1}\sum_{i=1}^n\sum_{t=1}^N b_{it}^k+O(\epsilon_k^2). \nonumber
\end{align}
where $\tilde \epsilon_k$ stands for $\epsilon_k^2\log \frac{1}{\epsilon_k}$.
Theorem \ref{gammaNpeq} is established. $\Box$

\section*{Acknowledgment}

The authors are grateful to Professor C.-S. Lin for  his helpful  discussions and comments on the article.


\begin{thebibliography}{99}
\bibitem{aly} J. J. Aly, Thermodynamics of a two-dimensional self-gravitating system,
Phy. Rev. A 49, No. 5, Part A (1994), 3771--3783.

\bibitem{barto1} D. Bartolucci, Compactness result for periodic multivortices in the electroweak theory,
Nonlinear Anal., 53 (2003), 277--297.

\bibitem{barto3} D. Bartolucci and G. Tarantello, Liouville type equations with singular data and their applications to periodic multivortices for the electroweak theory,
Comm. Math. Phys., 229 (2002), 3--47.

\bibitem{bennet} W. H. Bennet, Magnetically self-focusing streams,
Phys. Rev. 45 (1934), 890--897.

\bibitem{caffarelli} L. A. Caffarelli, Y. Yang, Vortex condensation in the Chern-Simons Higgs model: an existence theorem,
Comm. Math. Phys. 168 (1995), no. 2, 321--336.

\bibitem{lion1} E. Caglioti, P. L. Lions, C. Marchioro, C, M. Pulvirenti, A special class of stationary flows for two-dimensional Euler equations: a statistical mechanics description,
Comm. Math. Phys. 143 (1992), no. 3, 501--525.

\bibitem{lion2} E. Caglioti, P. L. Lions, C. Marchioro, C, M. Pulvirenti, A special class of stationary flows for two-dimensional Euler equations: a statistical mechanics description. II,
Comm. Math. Phys. 174 (1995), no. 2, 229--260.

\bibitem{ccl} S. A. Chang, C. C. Chen, C. S. Lin, Extremal functions for a mean field equation in two dimension. (English summary) Lectures on partial differential equations, 61–93,
New Stud. Adv. Math., 2, Int. Press, Somerville, MA, 2003.
\bibitem{chanillo1} S. Chanillo, M. K-H Kiessling, Rotational symmetry of solutions of some nonlinear problems in statistical mechanics and in geometry,
Comm. Math. Phys. 160 (1994), no. 2, 217--238.

\bibitem{chanillo2}S. Chanillo, M. K-H Kiessling, Conformally invariant systems of nonlinear PDE of Liouville type,
Geom. Funct. Anal. 5 (1995), no. 6, 924--947.

\bibitem{chchlin}J. L. Chern, Z. Y. Chen and C. S. Lin, Uniqueness of topological solutions and the structure of solutions for the Chern-Simons with two Higgs particles,
Comm. Math. Phys. 296 (2010), 323--351.

\bibitem{ChenLin1} C. C. Chen, C. S. Lin, Sharp estimates for solutions of multi-bubbles
in compact Riemann surfaces. Comm. Pure Appl. Math. 55 (2002), no. 6, 728--771.

\bibitem{chenlin2} C. C. Chen, C. S. Lin, Topological degree for a mean field equation on Riemann surfaces.
  Comm. Pure Appl. Math.  56  (2003),  no. 12, 1667--1727.


\bibitem{childress} S. Childress and J. K. Percus, Nonlinear aspects of Chemotaxis,
Math. Biosci. 56 (1981), 217--237.

\bibitem{CSW} M. Chipot, I. Shafrir, G. Wolansky, On the solutions of Liouville systems.  J. Differential Equations  140  (1997),  no. 1, 59--105.
\bibitem{CSW1} M. Chipot, I. Shafrir, G. Wolansky, Erratum: ``On the solutions of Liouville systems'' [J. Differential Equations 140 (1997), no. 1, 59--105; MR1473855 (98j:35053)].  J. Differential Equations  178  (2002),  no. 2, 630.

\bibitem{debye} P. Debye and E. Huckel, Zur Theorie der Electrolyte,
Phys. Zft 24 (1923), 305--325.

\bibitem{dunne} G. Dunne, Self-dual Chern-Simons Theories, Lecture Notes in Physics,
vol. m36, Berlin: Springer-Verlag, 1995.

\bibitem{dziar} J. Dziarmaga, Low energy dynamics of $[U(1)]^N$ Chern-Simons solitons and two dimensional nonlinear equations,
Phys. Rev. D 49 (1994), 5469--5479.

\bibitem{phys}L. Ferretti, S.B. Gudnason, K. Konishi,  Non-Abelian vortices and monopoles in SO(N) theories,
Nuclear Physics B, Volume 789, Issues 1-2, 21 January 2008, Pages 84-110.

\bibitem{gu-zhang} Y. Gu, L. Zhang, Degree counting theorems for singular Liouville systems, with Yi Gu.  to appear on Annali della Scuola Normale Superiore di Pisa, Classe di Scienze, DOI Number: $10.2422/2036-2145.201812_007$



\bibitem{huang-cvpde} H. Huang Existence of bubbling solutions for the Liouville system in a torus. Calc. Var. Partial Differential Equations 58 (2019), no. 3, Paper No. 99, 26 pp.

\bibitem{huang-zhang} H. Huang ; L. Zhang,
The domain geometry and the bubbling phenomenon of rank two gauge theory.
Comm. Math. Phys. 349 (2017), no. 1, 393–424.



\bibitem{HZ2021} H. Huang ; L. Zhang,
Local uniqueness of bubbling  solutions to the Liouville system, in preparation


\bibitem{hong} J. Hong, Y. Kim and P. Y. Pac, Multivortex solutions of the Abelian Chern-Simons-Higgs theory,
Phys. Rev. Letter 64 (1990), 2230--2233.

\bibitem{jackiw} R. Jackiw and E. J. Weinberg, Selfdual Chern Simons vortices,
Phys. Rev. Lett. 64 (1990), 2234--2237.


\bibitem{keller} E. F. Keller and L. A. Segel, Traveling bands of Chemotactic Bacteria: A theoretical analysis,
J. Theor. Biol. 30 (1971), 235--248.

\bibitem{kiess} M. K.-H. Kiessling, Statistical mechanics of classical particles with logarithmic interactions,
Comm. Pure Appl. Math. 46 (1993), no. 1, 27--56.

\bibitem{kiess2} M. K.-H. Kiessling, Symmetry results for finite-temperature, relativistic Thomas-Fermi equations,
Commun. Math. Phys. 226 (2002), no. 3, 607--626.

\bibitem{kiessling2}M. K.-H. Kiessling and J. L. Lebowitz, Dissipative stationary Plasmas: Kinetic Modeling Bennet Pinch, and generalizations,
Phys. Plasmas 1 (1994), 1841--1849.

\bibitem{kimleelee}C. Kim, C. Lee and B.-H. Lee, Schr\"odinger fields on the plane with $[U(1)]^N$ Chern-Simons interactions and generalized self-dual solitons,
Phys. Rev. D (3) 48 (1993), 1821--1840.

\bibitem{linzhang1} C. S. Lin, L. Zhang, Profile of bubbling solutions to a Liouville system.  Ann. Inst. H. Poincare Anal. Non Lineaire 27 (2010), no. 1, 117--143.

\bibitem{linzhang2} C. S. Lin, L. Zhang,   A topological degree counting for some Liouville systems of mean field equations, Comm. Pure Appl. Math. volume 64, Issue 4, pages 556--590, April 2011.

\bibitem{lin-zhang-jfa} C. S. Lin, L. Zhang, On Liouville systems at critical parameters, Part 1: One bubble. J. Funct. Anal. 264 (2013), no. 11, 2584--2636.

\bibitem{nolasco2} M. Nolasco, G. Tarantello, Vortex condensates for the SU(3) Chern-Simons theory,
Comm. Math. Phys. 213 (2000), no. 3, 599--639.

\bibitem{rubinstein} I. Rubinstein, Electro diffusion of Ions, SIAM, Stud. Appl. Math. 11 (1990).

\bibitem{spruck} J. Spruck, Y. Yang, Topological solutions in the self-dual Chern-Simons theory: existence and approximation,
Ann. Inst. H. Poincare Anal. Non Lineaire 12 (1995), no. 1, 75--97.

\bibitem{Wil} F. Wilczek,   Disassembling anyons.  Physical review letters 69.1 (1992): 132.

\bibitem{wolansky1} G. Wolansky, On steady distributions of self-attracting clusters under friction and fluctuations,
Arch. Rational Mech. Anal. 119 (1992), 355--391.

\bibitem{wolansky2} G. Wolansky, On the evolution of self-interacting clusters and applications to semi-linear equations with exponential nonlinearity,
J. Anal. Math. 59 (1992), 251--272.

\bibitem{wolansky3} G. Wolansky. Multi-components chemotactic system in the absence of conflicts. European Journal of Applied Mathematics, Volume 13, Issue 6, 2002.




\bibitem{yang} Y. Yang, Solitons in field theory and nonlinear analysis,
Springer-Verlag, 2001.

\bibitem{zhangcmp} L. Zhang, Blowup solutions of some nonlinear elliptic equations
involving exponential nonlinearities.
Comm. Math. Phys. 268 (2006), no. 1, 105--133.

\bibitem{zhangccm} L. Zhang, Asymptotic behavior of blowup solutions for elliptic equations with exponential nonlinearity and singular data.  Commun. Contemp. Math.  11  (2009),  no. 3, 395--411.

\end{thebibliography}
\end{document}